\documentclass[reqno]{amsart}

\usepackage[usenames,dvipsnames]{color}
\usepackage{hyperref}
\hypersetup{
    colorlinks=true,                            
    linkcolor=blue, 
    citecolor=blue,                             
    filecolor=blue,                                     
    urlcolor=blue                               
}

\usepackage{amsmath,amssymb,hyperref,amsthm}
\usepackage[all]{xy}

\newtheorem*{thm}{Theorem}
\newtheorem*{prop}{Proposition}
\newtheorem*{lem}{Lemma}
\newtheorem*{cor}{Corollary}

\newenvironment{pf}{\paragraph{{\sc Proof}}}{\qed\par\medskip}
\theoremstyle{definition}

\newtheorem*{rem}{Remark}
\newtheorem*{example}{Example}

\numberwithin{equation}{section}

\newcommand{\power}[2]{#1[\![#2]\!]}

\renewcommand{\ss}[1]{\scriptscriptstyle{#1}}
\newcommand{\op}{\ss{\mathrm{op}}}
\newcommand{\baff}{B^a}
\newcommand{\ebaff}{B^e}
\newcommand{\braid}[4]{\underbrace{#1_{#2}#1_{#3}#1_{#2}\cdots}
_{#4\text{ times}}}
\newcommand{\ex}[1]{#1^{\vee}}

\newcommand{\lp}{\left(}
\newcommand{\rp}{\right)}
\newcommand{\la}{\left\langle}
\newcommand{\ra}{\right\rangle}

\newcommand{\g}{\mathfrak{g}}
\newcommand{\h}{\mathfrak{h}}
\newcommand{\gl}{\mathfrak{gl}}
\newcommand{\Lsl}{\mathfrak{sl}}

\newcommand{\Sym}{\mathfrak{S}}

\newcommand{\sfv}{\mathsf{v}}
\newcommand{\sfw}{\mathsf{w}}

\newcommand{\E}{\mathcal{E}}
\newcommand{\F}{\mathcal{F}}
\newcommand{\K}{\mathcal{K}}

\newcommand{\V}{\mathcal{V}}
\newcommand{\HH}{\mathcal{H}}

\newcommand{\M}{\mathcal{M}}

\newcommand{\PP}{\mathcal{P}}

\newcommand{\W}{\mathcal W}

\newcommand{\C}{\mathbb{C}}
\newcommand{\nC}{\mathbb{C}^{\times}}
\newcommand{\N}{\mathbb{Z}_{\geq 0}}

\newcommand{\IP}{\mathbb{P}}

\newcommand{\Z}{\mathbb{Z}}
\newcommand{\qS}{\mathbb{S}}

\newcommand {\ol}[1]{\overline{#1}}
\newcommand {\ul}[1]{\underline{#1}}

\newcommand {\Ker}{\operatorname{Ker}}
\newcommand{\End}{\operatorname{End}}
\newcommand{\Hom}{\operatorname{Hom}}

\newcommand {\aand}{\qquad\text{and}\qquad}
\newcommand {\rk}{\operatorname{rk}}
\newcommand {\Aut}{\operatorname{Aut}}
\newcommand {\GL}{\mathrm{GL}}

\newcommand{\qW}{\lambda}
\newcommand{\cplim}{\mathrm{C}}

\newcommand{\LM}{\mathcal{L}} 
\newcommand{\IR}{I}    
\newcommand{\ws}{\mu}  

\newcommand {\fd}{finite--dimensional }

\newcommand {\wrt}{with respect to }

\newcommand{\ds}{\displaystyle}

\newcommand{\wt}[1]{\widetilde{#1}}
\newcommand{\wh}[1]{\widehat{#1}}

\newcommand {\comment}[1]{\footnote{\textcolor{blue}{#1}}}
\newcommand {\Omit}[1]{}

\newcommand{\aff}{\scriptscriptstyle{\operatorname{aff}}}

\newcommand {\ad}{\operatorname{ad}}
\newcommand {\Ad}{\operatorname{Ad}}

\newcommand {\qloop}{U_q(L\g)}
\newcommand {\qloopkm}{U_q^{\scriptscriptstyle{\mathrm{KM}}}(L\g)}
\newcommand {\qlsl}[1]{U_{q_{#1}}^{\scriptscriptstyle{\mathrm{KM}}}(L\Lsl_2)}

\newcommand {\vembed}{\varphi}
\newcommand {\bsign}{o}        
\newcommand {\biso}[1]{\beta_{#1}} 
\newcommand {\cg}{\psi} 
\newcommand {\cgn}{\ol{\psi}} 
\newcommand {\GK}{\gamma}

\newcommand {\qloopsl}[1]{U_q(L\Lsl_{#1})}

\newcommand{\Rloop}{\operatorname{Rep}_{\scriptscriptstyle{\text{fd}}}(\qloop)}
\newcommand{\Rloopkm}{\operatorname{Rep}_{\scriptscriptstyle{\text{fd}}}(\qloopkm)}

\newcommand {\isom}{\stackrel{\sim}{\rightarrow}}

\newcommand{\cbin}[2]{\left(\begin{array}{c} #1\\ #2\end{array}\right)}
\newcommand {\bin}[3]{\left[\begin{array}{c}#1\\#2\end{array}\right]_{#3}}
\newcommand {\qbin}[3]{\left[\begin{array}{c}#1\\#2\end{array}\right]_{#3}}

\newcommand {\eg}{{\it e.g.,} }
\newcommand {\ie}{{\it i.e.,} }

\newcommand {\bfA}{{\mathbf A}}
\newcommand {\bfI}{{\mathbf I}}
\newcommand {\bfE}{{\mathbf E}}
\newcommand {\bfJ}{{\mathbf J}}

\renewcommand {\H}{\mathcal H}

\newcommand {\Id}{\operatorname{Id}}
\renewcommand {\det}{\operatorname{det}}

\renewcommand {\gl}{\mathfrak{gl}}
\renewcommand {\sl}{\mathfrak{sl}}

\newcommand{\Nqv}{\Phi(\sfw)}
\newcommand {\sj}{_{s}}
\newcommand {\uj}{_{u}}
\newcommand {\qWg}{quantum Weyl group }
\newcommand {\qloopo}{U_q^0(L\g)}
\newcommand {\lambdao}{\qW_{\bsign}}

\newcommand {\nofootnote}[1]{}
\newcommand {\nocomment}[1]{}

\title
[Quantum Weyl group action of the coroot lattice]
{An abelian formula for the quantum Weyl group action of the coroot lattice}
\author[S. Gautam]{Sachin Gautam}
\address{Department of Mathematics,
The Ohio State University,
231 W 18th Ave,
Columbus, OH 43210}
\email{gautam.42@osu.edu}
\author[V. Toledano Laredo]{Valerio Toledano Laredo}
\address{Department of Mathematics,
Northeastern University,
567 Lake Hall,
360 Huntington Avenue,
Boston, MA 02115}
\email{V.ToledanoLaredo@neu.edu}

\begin{document}
\maketitle

\begin{abstract}
Let $\g$ be a complex simple Lie algebra and $\qloop$ its quantum loop 
algebra, where $q\in\nC$ is not a root of unity. We give an explicit formula for
the quantum Weyl group action of the coroot lattice $Q^\vee$ of $\g$ on \fd
representations of $\qloop$ in terms of its commuting generators. The answer
is expressed in terms of the Chari--Pressley series, whose evaluation on
highest weight vectors gives rise to Drinfeld polynomials. It hinges on a strong
rationality result for that series, which is derived in the present paper. As an
application, we identify the action of $Q^\vee$ on the equivariant $K$--theory
of Nakajima quiver varieties with that of explicitly given determinant line
bundles.
\end{abstract}

\Omit{arXiv abstract:
Let g be a complex simple Lie algebra and Uq(Lg) its quantum loop algebra,
where q is not a root of unity. We give an explicit formula for the quantum Weyl
group action of the coroot lattice Q of g on finite-dimensional representations
of Uq(Lg) in terms of its commuting generators. The answer is expressed in
terms of the Chari-Pressley series, whose evaluation on highest weight vectors
gives rise to Drinfeld polynomials. It hinges on a strong rationality result for that
series, which is derived in the present paper. As an application, we identify the
action of Q on the equivariant K-theory of Nakajima quiver varieties with that
of explicitly given determinant line bundles.
}

\setcounter{tocdepth}{1}
\tableofcontents

\section{Introduction}\label{sec: intro}

\subsection{}

Let $\g$ be a complex simple Lie algebra, $q\in\nC$ an element of infinite order,
and $\qloop$ the quantum loop algebra of $\g$. If $V$ is a \fd 
representation
of $\qloop$, the quantum Weyl group action of the affine braid group $\baff$ on
$V$ restricts to one of the coroot lattice $Q^{\vee}\subset\baff$ which commutes
with the abelian generators $\{\cg_{i,\pm k}^{\pm}\}_{k\geq 0}$ of $\qloop$. Since 
these generate a maximal commutative subalgebra $\qloopo$ of $\qloop$, this
suggests that $Q^\vee$ might act on $V$ through elements of $\qloopo$. 
This is not clear a priori, however. For example, if $\g=\sl_2$, the coroot lattice
is generated by the element $S_0S_1\in\baff$, and its action on $V$ is given by
the product of six $q$--exponentials
\[\begin{split}
\qS_0 \qS_1 = &\exp_{q^{-1}}(q^{-1}\E_0\K_0^{-1})
\exp_{q^{-1}}(-\F_0)\exp_{q^{-1}}(q\E_0\K_0)\,
q^{\HH_0(\HH_0+1)/2}\\
\cdot&\exp_{q^{-1}}(q^{-1}\E_1\K_1^{-1})
\exp_{q^{-1}}(-\F_1)\exp_{q^{-1}}(q\E_1\K_1)\,
q^{\HH_1(\HH_1+1)/2}
\end{split}\]
which isn't readily expressed in terms of elements of $\qloopo$.

\subsection{} 

The goal of this paper is to obtain explicit abelian formulae for the \qWg action of
$Q^\vee$. When $q$ is formal and $\g=\Lsl_2$ or $\gl_2$, a formula expressing
this action in terms of the primitive--like loop generators $\{H_{i,k}\}$ was derived in
our earlier work \cite{sachin-valerio-3}, and used to prove the monodromy
conjecture of the second author \cite{valerio3}. In the present work, we obtain a
formula in terms of the group--like generators $\{\psi_{i,k}^\pm\}$ which is new even
for $\g=\sl_2$. Interestingly, the answer is expressed as the evaluation of a well--known
and very basic object, namely the series introduced by Chari--Pressley, whose
truncation on highest weight vectors gives rise to Drinfeld polynomials. When
$\g$ is simply--laced and $V$ is the equivariant $K$--theory of a Nakajima
quiver variety, we rely on our new formulae to describe the action of $Q^\vee$
on $V$ as tensoring with explicit determinant line bundles.

\subsection{}

To state our results more precisely, recall that $\qloop$ admits two distinct
presentations, which were shown to be isomorphic by Beck \cite{beck-braid}.
To distinguish them, we denote by $\qloopkm$ the algebra given in the Kac--Moody
presentation, with generators $\{\E_i,\F_i,\K_i^{\pm 1}\}_{i\in\wh{\bfI}}$, where
$\bfI$ is the set of vertices of the Dynkin diagram of $\g$ and $\wh{\bfI}=\bfI
\sqcup\{0\}$, and by $\qloop$ the algebra presented on the loop generators
$\{E_{i,k},F_{i,k},\cg^{\pm}_{i,\pm\ell}\}_{i\in\bfI,k\in\Z,\ell\in\Z_{\geq 0}}$ \cite
{drinfeld-yangian-qaffine}.

We denote Beck's isomorphism by $\biso{\bsign}:\qloop\to\qloopkm$. It depends
on a choice of a sign $\bsign:\bfI\to\{\pm 1\}$ which is bipartite in the sense that
$\bsign(i)\bsign(j)=-1$ if $i,j\in\bfI$ are connected in the Dynkin diagram of $\g$.
Under this isomorphism, $\cg_{i,0}^{\pm}$ are identified with $\K_i^{\pm 1}$, for
all $i\in\bfI$.

\subsection{}\label{issec:QW} 

Let $V$ be a finite--dimensional type I representation of $\qloop$, and view it as
a representation of $\qloopkm$ via $\biso{\bsign}$. For any $i\in\wh{\bfI}$, let
$\HH_i\in\End(V)$ be the semisimple operator with integer eigenvalues such
that $\K_i=q_i^{\HH_i}$, where $q_i=q^{d_i}$ and $\{d_i\}_{i\in\wh{\bfI}}$ are
the symmetrizing integers. The triple $q$--exponentials
\[
\qS_i = \exp_{q_i^{-1}}(q_i^{-1}\E_i\K_i^{-1})
\exp_{q_i^{-1}}(-\F_i)\exp_{q_i^{-1}}(q_i\E_i\K_i)
q_i^{\HH_i(\HH_i+1)/2}
\]
define an action of $\baff$ on $V$ \cite{kirillov-reshetikhin,lusztig-canonicalII,
lusztig-book,saito-PBW,soibelman}. Here,
\[\exp_q(x)=\sum_{n=0}^{\infty} q^{n(n-1)/2}\frac{x^n}{[n]!},
\qquad
[n]!=[n]\cdots [1]
\aand
[k]=\frac{q^k-q^{-k}}{q-q^{-1}}\]
We denote this action by $\lambdao : \baff\to\GL(V)$. 

Recall that the coroot lattice $Q^{\vee}$ is a subgroup of $\baff$
\cite[Ch.~3]{macdonald-affine}, and let $L_i\in\baff$ be the element
corresponding to the coroot $\alpha_i^{\vee}$. The main goal of
this paper is to give a formula for $\lambdao(L_i)\in\GL(V)$ in
terms of the commuting generators $\{\cg^{\pm}_{i,\pm k}\}
_{i\in\bfI,k\in\Z_{\geq 0}}$. 

\subsection{The Chari--Pressley series}\label{issec:q-diff}

For any $i\in\bfI$, set
\[\cg_i^{\pm}(z) = \sum_{n=0}^{\infty} \cg^{\pm}_{i,\pm n} z^{\mp n}
\in \power{U_q^0(L\g)}{z^{\mp 1}}
\aand
\cgn_i^{\pm}(z) = \cg_{i,0}^{\mp}\cdot \cg_i^{\pm}(z)\]
so that $\cgn_i^+(\infty)=1=\cgn_i^-(0)$. Let $\PP_i^{\pm}(z)\in 1+z^
{\mp 1}\power{U_q^0(L\g)}{z^{\mp 1}}$ be the unique formal solution
of the $q$--difference equation

\begin{equation}\label{ieq:diff}
\PP_i^{\pm}(q_i^2z) = \cgn_i^{\pm}(z)\PP_i^\pm(z)
\end{equation}

The series $\PP_i^\pm(z)$ were introduced by Chari--Pressley \cite[\S 3.5]
{cp-qaffine}, and shown to truncate on the highest weight space of a \fd 
representation $V$ of $\qloop$, thus obtaining the existence of Drinfeld
polynomials.\footnote{Our conventions differ slightly from those of {\em loc.
cit.} but agree, up to changing $z$ by $z^{-1}$ for $\PP^+(z)$, with the ones
from \cite[\S 3]{cp-rootofunity}.} Its classical analogue and the corresponding
truncation result were obtained earlier in Chari's work on integrable representation
of affine Lie algebras \cite[Prop.~1.1]{chari-integrable}.

On an arbitrary weight space of $V$, the form of the eigenvalues of $\cg_i^
\pm(z)$ given by Frenkel--Reshetikhin \cite[Prop.~1]{frenkel-reshetikhin-qchar}
implies, and is in fact easily seen to be equivalent to, the rationality of the
eigenvalues of $\PP_i^\pm(z)$ (see Section \ref{issec:evs}). 

A stronger result holds when $\g$ is simply--laced and $V$ is the equivariant
$K$--theory of a Nakajima quiver variety. In that case, the action of $\PP_i^+(z),
\PP_i^-(z)$ given in \cite{nakajima-qaffine} shows that they are the Taylor
expansions at $z=\infty,0$ of rational functions (see \ref{ssec:U0-action}--\ref
{ssec:nakajima-lattice}).

\subsection{Main result}\label{issec:mainthm} 

In this paper, we prove the rationality of $\PP_i^\pm(z)$ for an arbitrary $\g$ and
$V$. We then express $\lambdao(L_i)$ as the normalised limit of $\PP_i^+(z)$
as $z\to 0$ (Theorem \ref{thm:cp-series}).\Omit{\footnote{An {\em a posteriori}
explanation of this can be inferred from their commutation relation with the raising
/lowering operators (see Lemma \ref{lem:comm} and Section \ref{ssec:pf-3-comm}).}}

\Omit{
\subsection{}\label{issec:mainthm} 

Our main result (Theorem \ref{thm:cp-series}) hinges on, and establishes
the rationality of $\PP_i^\pm(z)$ for an arbitrary $\g$ and $V$. It then
expresses $\lambdao(L_i)$ as a normalised limit of $\PP_i^+(z)$ as
$z\to 0$.\Omit{\footnote{An {\em a posteriori} explanation of this can
be inferred from their commutation relation with the raising/lowering
operators (see Lemma \ref{lem:comm} and Section \ref{ssec:pf-3-comm}).}}
}

\begin{thm}\label{ithm:main}
Let $V$ be a finite--dimensional type I representation of $\qloop$.
\begin{enumerate}
\item\label{iit:rat} The action of $\PP_i^+(z)$ (resp. $\PP_i^-(z)$) on
$V$ is the Taylor series at $z=\infty$ (resp. $z=0$) of a rational
$\End(V)$--valued function $P_i^+(z)$ (resp. $P_i^-(z)$).
\item\label{iit:lim} There is an element $\cplim_i\in\GL(V)$ which commutes
with the action of $\qloopo$, and such that $\ds z^{\H_i}P_i^+(z)=\cplim_i
P_i^-(z)$. It follows that
\begin{equation}\label{eq:Ci}
\cplim_i = \lim_{z\to 0} z^{\H_i}P_i^+(z)
= \lim_{z\to\infty} z^{\H_i}P_i^-(z)^{-1}
\end{equation}

\item\label{iit:lattice} 
The quantum Weyl group action of the generator $L_i\in Q^\vee$ on $V$
is given by
\[\lambdao(L_i) = \left(\bsign(i)q_i\right)^{\H_i}\cdot \cplim_i^{-1}\]
\end{enumerate}
\end{thm}

We give several applications of Theorem \ref{ithm:main} in \ref{issec:evs}--\ref
{ss:CKL} below, and sketch its proof in \ref{ss:sketch start}--\ref{ss:sketch end}.

\subsection{Remarks}

\begin{enumerate}
\item
Parts (1)--(2) of Theorem \ref{ithm:main} may be interpreted from the point
of view of $q$--difference equations as follows. If $V$ is representation of $\qloop$
with \fd weight spaces, $\cgn^{\pm}_i(z)$ are rational $\End(V)$--valued functions
\cite{beck-kac,hernandez-drinfeld-coproduct,sachin-valerio-2}. If $|q|\neq 1$,
the unique solutions $P_i^{\pm}(z)$ of \eqref{ieq:diff} define holomorphic $\GL
(V)$--valued functions in a neighbourhood of $z=\infty/0$. These possess a
meromorphic continuation to $\nC$, but may have an essential singularity at
$z=0/\infty$.

If $V$ is finite--dimensional, however, Theorem \ref{ithm:main} \eqref{iit:lim} rules
out the presence of essential singularities and implies that $P_i^{\pm}(z)$ have
at worst a pole at $z=0/\infty$, whose order is the eigenvalue of $\HH_i$. In
particular, the $q$--difference equation \eqref{ieq:diff} has trivial monodromy
on $V$.

\item In \cite[(5.31)]{frenkel-hernandez-recent}, Frenkel--Hernandez recently showed
that the series $\PP_i^-(z)$ can be expressed as a Laurent monomial in a family of
formal power series $\{X_j(z)\}_{j\in\bfI}$ they introduced earlier \cite[Prop.~5.5]{frenkel-hernandez}.\nofootnote
{The series $T_i(z)$ is denoted by $X_i(z)$ in \cite{frenkel-hernandez-recent} to distinguish
it from quantum Weyl group operators.} They also show that on a (not necessarily
finite--dimensional) highest weight representation, the action of an appropriate normalisation
of $X_j(z)$ is {\it polynomial} \cite[Thm.~5.9]{frenkel-hernandez}.

While it is of a similar nature, this result is essentially logically independent of Theorem
\ref{ithm:main}. Indeed, the rationality of $\PP_i^{\pm}(z)$ does not hold for an arbitrary
category $\mathcal{O}$ representation.{\footnote{Combining \cite[Thm.~5.5]{frenkel-hernandez}
and \cite[(5.31)]{frenkel-hernandez-recent}, one can show that the {\em normalised}
Chari--Pressley series are rational on highest weight representations, which is not
true without the normalisation.}} Conversely, the rationality results of \cite{frenkel-hernandez,
frenkel-hernandez-recent} 
do not apply to general finite--dimensional representations, and cannot therefore
be deduced from Theorem \ref{ithm:main}.

Finally, if $V$ is a representation which is both highest weight and finite--dimensional,
Theorem \ref{ithm:main} \eqref{iit:rat} can be proved very easily as follows. 
Applying the classification of finite--dimensional irreducible representations
in terms of Drinfeld polynomials to the socle of $V$, one obtains the rationality
of $\PP_i^{\pm}(z)$ on its highest weight space. A simple induction combined
with the commutation relations between $\PP^{\pm}_i(z)$ and lowering operators
(see Lemma \ref{lem:comm} below) then shows that the same on the whole of $V$.\\

\end{enumerate}

\subsection{Eigenvalues of $\cg_i^\pm(z)$}\label{issec:evs}

As a first application of Theorem \ref{ithm:main}, we give an alternative determination
of the eigenvalues of $\cg_i^\pm(z)$ on a type I \fd representation $V$ of $\qloop$
obtained in \cite[Prop.~1]{frenkel-reshetikhin-qchar}, which is independent of the
classification of irreducible ones. We also compute the corresponding eigenvalues
of the operators $\lambdao(L_i)$.

\begin{prop}\label{ipr:evs}\hfill 
\begin{enumerate}
\item The eigenvalues of $\cg_i^\pm(z)$ on $V$ are of the form
\Omit{
\begin{equation}\label{eq:}
q_i^{-\deg Q_i+\deg R_i}\frac{Q_i(q_i^2z)R_(z)}{Q_i(z)R_i(q_i^2 z)}
\end{equation}
where $Q_i,R_i$ are monic polynomials in $z$ with roots in $\C^\times$.
}
\begin{equation}\label{eq:eval psi}
q_i^{-\deg r_i}\left(\frac{r_i(q_i^2 z)}{r_i(z)}\right)^\pm
\end{equation}
where $r_i(z)$ is a rational function such that $r_i(0)\in\C^\times$, and $(-)^
\pm$ denotes the Taylor expansion at $z=\infty,0$ respectively.
\item Normalise $r_i(z)$ by requiring that $\lim_{z\to\infty}r_i(z)z^{-\deg r_i}=1$.
Then, on the corresponding generalised eigenspace of $\cg_i^\pm(z)$,
$\lambdao(L_i)$ has a single eigenvalue given by $(\bsign(i)q_i)^{\H_i}/r_i(0)$.
Explicitly, if 
\[r_i(z)=\prod_{j=1}^{m} (z-a_{i,j})/\prod_{k=1}^n (z-b_{i,k})\]
where $a_{i,j},b_{i,k}\in\C^\times$, the corresponding eigenvalue of
$\lambdao(L_i)$ is given by 
\[(-\bsign(i)q_i)^{m-n}\prod_k b_{i,k}/\prod_j a_{i,j}\]
\end{enumerate}
\end{prop}
\begin{pf}
(1) Let $\cg_i:\C\to GL(V)$ the rational function whose Taylor expansions at
$z=\infty,0$ are equal to $\cg^+(z),\cg^-(z)$ respectively \cite{beck-kac,
hernandez-drinfeld-coproduct,sachin-valerio-2}. Let $P_i^+(z)=P_i^+(z)
\sj\cdot P_i^+(z)\uj$ be the Jordan decomposition of $P_i^+(z)$. Since
$[P_i^+(z),P_i^+(z')]=0$, it follows from \eqref{ieq:diff} that the Jordan
decomposition of $\cg_i(z)$ is given by
\[\cg_i(z)\sj = q_i^{+\HH_i}P_i^{+}(q_i^2z)\sj P_i^+(z)\sj^{-1}
\aand
\cg_i(z)\uj = P_i^{+}(q_i^2z)\uj P_i^+(z)\uj^{-1}\]
The commutativity of $P_i^+(z)$ and the fact that it is rational by Theorem
\ref{ithm:main} imply the same for its semisimple and unipotent components
(see, \eg \cite[Lemma 4.12]{sachin-valerio-2}). In particular,
its eigenvalues are rational functions of $z$
taking the value $1$ at $z=\infty$.

Let $p_i(z)$ be an eigenvalue of $P_i^+(z)$, and define $r_i(z)\in\C(z)$
by $p_i(z)=r_i(z)z^k$ and $r_i(0)\neq 0$. Then, $\psi_i(z)\sj$ acts on the
corresponding generalised eigenspace $V[p_i]$ as multiplication by
\[ q_i^{\HH_i}
q_i^{-2\deg r_i}\frac{r_i(q_i^2 z)}{r_i(z)} \]
Since $\psi(\infty)=\psi(0)^{-1}$, it follows that $q_i^{\HH_i}=q_i^{\deg r_i}$,
as claimed.

(2) Let $C_i\in GL(V)$ be given by \eqref{eq:Ci}. Then, $C_i$ has a single
eigenvalue on $V[p_i]$ given by $r_i(0)$, and the claim follows by Theorem
\ref{ithm:main} \eqref{iit:lattice}.
\end{pf}

\subsection{Relation with the Drinfeld coproduct}

In \cite{drinfeld-yangian-qaffine}, Drinfeld introduced a topological coproduct
$\Delta_D$ on $\qloop$. Composing it with $\tau_z\otimes\Id$, where
$\tau_z:\qloop\to\qloop[z^{\pm 1}]$ is the shift homomorphism, 
yields an algebra homomorphism $\Delta_{D,z}:\qloop
\to\qloop^{\otimes 2}((z))$ called the {\em deformed} Drinfeld coproduct, which
was introduced and studied by Hernandez \cite[\S 6.2]{hernandez-05}.

Given two finite--dimensional representations $V_1$ and $V_2$, the action of
$\qloop$ on $V_1\otimes V_2$ via $\Delta_{D,\zeta}$ is rational in $\zeta$ \cite
[Lemma~3.20]{hernandez-drinfeld-coproduct} (see also \cite[Thm.~4.3]{GTL17}).
We denote this representation by $V_1\otimes_\zeta V_2$. It is defined for all
but finitely many $\zeta\in\nC$, the finite set being dependent on $V_1$ and
$V_2$ \cite[Thm.~4.3 (iii)]{GTL17}.

The following result states roughly that the lattice operators are group--like \wrt
the deformed Drinfeld coproduct.

\begin{prop}
Let $V_1,V_2\in\Rloop$ be of type $I$. Then,
\[
\qW_{V_1\otimes_\zeta V_2,\bsign}(L_i) = \zeta^{-\HH_i}\cdot
\qW_{V_1,\bsign}(L_i)\otimes \qW_{V_2,\bsign}(L_i)
\]
\end{prop}

\begin{pf}
Recall that 
\[\Delta_{D,\zeta}(\psi_i^{\pm}(z)) = \psi_i^{\pm}(\zeta^{-1}z)\otimes
\psi_i^{\pm}(z)\]
Equation \eqref{ieq:diff} then implies that
\[
\Delta_{D,\zeta}(\PP_i^{\pm}(z)) = \PP_i^{\pm}(\zeta^{-1}z)\otimes \PP_i^{\pm}(z)
\]
Thus, the normalised limit $\cplim_{i,V_1\otimes_\zeta V_2}$
from \eqref{eq:Ci} is also group--like, as computed below.
\[
\begin{aligned}
\cplim_{i,V_1\otimes_\zeta V_2} &= \lim_{z\to 0} z^{\HH_i}\otimes z^{\HH_i}
\ \circ\ 
\PP_i^+(\zeta^{-1}z) \otimes \PP_i^+(z) \\
&= \zeta^{\HH_i}\otimes \Id \ \circ\ \lim_{z\to 0}
(\zeta^{-1}z)^{\HH_i}\PP_i^+(z) \otimes z^{\HH_i}\PP_i^+(z)\\
&= \zeta^{\HH_i}\otimes \Id\ \circ\ \cplim_{i,V_1}\otimes\cplim_{i,V_2}
\end{aligned}
\]
The result now follows from Theorem \ref{ithm:main} \eqref{iit:lattice}.
\end{pf}

\subsection{An explicit formula}

We next give an explicit formula for the operators $\lambdao(L_i)$ in terms
of the operators $\{H_{i,k}\}_{i\in\bfI,k\in\Z}$. Here, for a fixed $i\in\bfI$,
$\{H_{i,k}\}_{k\neq 0}$ are defined by
\[\psi_i^{\pm}(z) = \psi_{i,0}^{\pm} \exp\lp\pm (q_i-q_i^{-1})\sum_{r=1}^{\infty}
H_{i,\pm r} z^{\mp r}\rp\]
and $H_{i,0} = \HH_i$ is the unique semisimple operator with $\Z$--eigenvalues
such that $\psi_{i,0}^{\pm} = \K_i^{\pm 1} = q_i^{H_{i,0}}$ as in Section
\ref{issec:QW} above. Note that, with the conventions followed in this paper,
$\HH_i$ or $H_{i,0}$ is not an element of $U_q(L\g)$, but only a well--defined
operator on finite--dimensional, type I representations.

For every $r\geq 1$, set $\ds
\wt{H}_{i,r} = H_{i,0}+\sum_{s=1}^{r} (-1)^s \cbin{r}{s}
\frac{s}{[s]_i} H_{i,s}$.
The following change of variables is easy to verify (see Section \ref{ssec:euler})
\[
\left.(1-q_iz^{-1})^{\HH_i}\PP_i^+(z)^{-1}\right|_{z=q\frac{t-1}{t}} 
= \exp\lp\sum_{r=1}^{\infty} \wt{H}_{i,r} \frac{t^r}{r}\rp
\]
Theorem \ref{ithm:main} \eqref{iit:lattice} then implies the following result (see
Corollary \ref{cor:euler})

\begin{prop}
Let $V$ be a finite--dimensional, type I representation of $\qloop$. Then,
\[
\lambdao(L_i) = (-\bsign(i))^{H_{i,0}} \cdot \lim_{t\to 1}
\exp\lp \sum_{r=1}^{\infty} \wt{H}_{i,r} \frac{t^r}{r}\rp
\]
\end{prop}

In the formal $\hbar$--adic setting, and when $\g=\sl_{2}$, this formula was
obtained in our earlier work \cite{sachin-valerio-3}. Note that the classical
limit of $\wt{H}_{i,r}$ is the element of $\g[u]$ given by
\[
\wt{h}_{i,r} = h_i\otimes \lp\sum_{s=0}^r (-1)^s \cbin{r}{s} u^s\rp
= h_i\otimes (1-u)^r
\]
Thus, the expression written above is equal to
\[
\lim_{t\to 1} \exp\lp h_i\otimes \sum_{n=1}^{\infty} \frac{(1-u)^nt^n}{n}\rp
= \lim_{t\to 1} \exp\lp -h_i\otimes \log(1-t(1-u))\rp
\]
which is the formal expansion of $u^{-h_i}$ in the loop group of $\g$. It agrees
with the image of $\alpha_i^{\vee}$ via the embedding of (the Tits extension of)
the affine Weyl group into the corresponding loop group, as can be verified by
an easy $2\times 2$ matrix calculation.

\subsection{Quiver varieties}

As another application of Theorem \ref{ithm:main}, we compute the \qWg action
of $Q^\vee$ on the equivariant $K$--theory of Nakajima quiver varieties, assuming
$\g$ is simply--laced. Specifically, 
we compose Beck's
isomorphism with Nakajima's action of $\qloop$ on these spaces \cite{nakajima-qaffine}
\[\qloopkm \stackrel{\biso{\bsign}^{-1}}{\longrightarrow}
\qloop \stackrel{\Nqv}{\longrightarrow}
\End\lp\bigoplus_{\sfv} K_{\GL(\sfw)\times\nC}(\M(\sfv,\sfw))
\rp\]
The restriction of $\Nqv$ to 
$U_q^0(L\g)$ is given in terms of the following complex
of tautological vector bundles on $\M(\sfv,\sfw)$ (see Section \ref{ssec:bundles}
for details)
\[
C_k(\sfv,\sfw):\quad
q^{-2}\V_k \longrightarrow
q^{-1}\lp \W_k \oplus\bigoplus_{\ell:a_{k\ell}=-1} \V_\ell\rp
\longrightarrow
\V_k
\]
where $(a_{k\ell})$ is the Cartan matrix of $\g$. We then prove that the following
holds (see Theorem \ref{thm:nakajima-lattice})
\begin{thm}\label{ithm:qv}
If $\g$ is simply--laced, the action of $\bsign(k)^{\H_k}L_k$ on $K_{\GL(\sfw)\times
\nC}(\M(\sfv,\sfw))$ is given by tensoring with the line bundle
\begin{equation}\label{eq:det line}
\det(C_k(\sfv,\sfw))^* = 
q^{\rk(C_k(\sfv,\sfw))}
\det(\W_k)^*\otimes\det(\V_k)^{\otimes 2}
\otimes\lp\bigotimes_{\ell:a_{k\ell}=-1} \det(\V_\ell)^*\rp
\end{equation}
\end{thm}

\subsection{Comparison with work of Cautis--Kamnitzer--Licata}\label{ss:CKL}
In \cite[\S 7]{CKL}, the authors give an action of the extended affine braid
group $\ebaff$ on the {\it non--equivariant} derived category of coherent
sheaves on Nakajima quiver varieties. This is achieved in two steps

\begin{enumerate}
\item The categorical action of $U_q(\g)$ is obtained via derived functors
coming from Nakajima's Hecke correspondences. This yields in particular
an action of the finite braid group through twisting with Rickard complexes
\`{a} la Chuang--Rouquier \cite[Thm.~6.4]{chuang-rouquier}. See \cite[Thm.~7.1]{CKL}.
\item For the coweight lattice $P^\vee\subset\ebaff$, the action of the generator
$Y_i$ corresponding to the fundamental coweight $\varpi_i^{\vee}$ 
is given by tensoring with $\det(\V_i)$. 
\end{enumerate}

As pointed out in \cite[Thm.~7.3]{CKL}, these complexes and line bundles do
{\it not} give rise to an action of $\ebaff$ in the equivariant setting since one of
its defining relations, namely
\[T_i^{-1}Y_iT^{-1}=Y_i^{-1}\prod_{j\neq i} Y_j^{-a_{ij}}\]
only holds up to shifts in the equivariant derived category. In particular, one
does not get a categorification of the \qWg action of the affine braid group
$\ebaff$ on the equivariant $K$--theory.

Theorem \ref{ithm:qv} and the fact that $L_k = Y_k^2\prod_{\ell:a_{k\ell}
=-1} Y_\ell^{-1}$ imply that the line bundles used in \cite{CKL} differ from those
given by Theorem \ref{ithm:qv} by
tensoring with $q^{\rk(C_k(\sfv,\sfw))}\det(\W_k)^*$ and correcting by a bipartite
sign. We conjecture that Rickard complexes and $\det(C_k(\sfv,\sfw))^*$ give rise to an action of $\baff$ on the equivariant
derived category of quiver varieties. This would then give a categorification
of the \qWg action of $\baff$ by Theorem \ref{ithm:qv}.

\subsection{Extended affine braid group actions}

It is a natural to ask whether the quantum Weyl group action of $\baff$ on a \fd 
representation $V$ of $\qloop$ canonically extends to $\ebaff$. The answer
appears to be no, unless one makes consistent choices of $d$th roots where
$d$ is the exponent of the quotient $P^{\vee}/Q^{\vee}$.

For instance, if $V$ is irreducible, with Drinfeld polynomials $\{Q_i(z)\}_{i\in\bfI}$
then, by Proposition \ref{ipr:evs} and up to a sign and a power of $q$, the coroot
lattice operator $L_i$ acts on the highest weight subspace as multiplication by $Q_
i(0)^{-1}$. Defining the action of the coweight lattice operator $Y_i= \lp\prod_j L_j
^{c_{ij}}\rp^{1/d}$ on that subspace therefore involves choosing a $d$th roots of 
each $\prod_j Q_j(0)^{c_{ij}}$, $i\in\bfI$, where $(c_{ij})_{i,j\in\bfI}$ is the integer
matrix obtained from the relations $d\varpi_i^{\vee} = \sum_{j\in\bfI} c_{ij}\alpha_
j^{\vee}$. In the geometric context, this amounts to introducing $d$th roots of the
line bundles $\bigotimes_j\det(\W_j)^{\otimes c_{ij}}$.

In \cite[\S 9.1]{etingof-varchenko}, this issue is addressed as follows. Given a \fd
representation $V$, $\qloop$ acts on $V[z^{\pm 1}]=V\otimes_\C\C[z^{\pm 1}]$
via the shift homomorphism $\tau_z : \qloop\to\qloop[z^{\pm 1}]$. Then, there is
a natural action of $\ebaff$ on $V[z^{\pm\frac{1}{d}}]$
extending that of $\baff$ on $V[z^{\pm 1}]$.

\subsection{Idea of the proof of Theorem \ref{ithm:main}}\label{ss:sketch start}
The proof of Theorem \ref{ithm:main} is based on a rank $1$ reduction and
a straightening identity for $\qloopsl{2}$, which might be of independent
interest.
The rank $1$ reduction relies on
the defining relations of $\ebaff$ (Corollary \ref{cor: rank1-formula})
and Beck's result \cite[Prop.~3.8]{beck-braid}, and is obtained
in Proposition \ref{pr:rk1-red}.

To motivate the straightening identity (Prop. \ref{pr:straight}), 
consider the classical and non--affine cases first.
If $\{e,f,h\}$ is the standard basis
of $\Lsl_2$, it is a straightforward exercise 
(\cite[Problem 2.15.1 (b)]{pasha-book}, or
\cite[Ex.~3.2]{kac}{\footnote{there is typo in the latter reference,
$h-m-n-2j$ there, should be $h-m-n+2j$.}})
to show that
\begin{equation}\label{ieq:pbw-sl}
e^{(m)}f^{(n)} = \sum_{p=0}^{\min(m,n)} \cbin{h-m+n}{p}
f^{(n-p)}e^{(m-p)}
\end{equation}
where, $\ds x^{(r)} = \frac{x^r}{r!}$
and $\ds\cbin{x}{r} = \frac{x(x-1)\cdots (x-r+1)}{r!}$ is
viewed as a polynomial in the variable $x$.
This implies that the following
relation holds in $\power{U(\Lsl_2)}{z^{-1}}$

\begin{equation}\label{ieq:series-sl}
\sum_{n=0}^{\infty} (-1)^n e^{(n)}f^{(n)}z^{-n} = 
p(z) \sum_{n=0}^{\infty} (-1)^n f^{(n)}e^{(n)}z^{-n}
\end{equation}
where
\[
p(z) = \sum_{\ell=0}^\infty (-1)^{\ell} \cbin{h}{\ell} z^{-\ell}
= (1-z^{-1})^h = \exp\lp -h\sum_{n=1}^{\infty} \frac{z^{-n}}{n}\rp
\]

Adapting this identity to $U_q(\Lsl_2)$ is direct: one merely has to replace $n$ by
$[n]=\frac{q^n-q^{-n}}{q-q^{-1}}$ everywhere. Thus, for instance, $[h-m] = \frac{q^
{-m}K - q^mK^{-1}}{q-q^{-1}}$.

\subsection{}
The analogue of \eqref{ieq:pbw-sl} for the loop algebra $U(L\Lsl_2)$ was obtained
by Garland \cite[Lemma~7.5]{garland}. Denote the standard basis vectors of
$\Lsl_2[t,t^{-1}]$ by $\{e_k,f_k,h_k\}_{k\in\Z}$. Then, the following 
holds{\footnote{This equation
differs slightly from the one in {\em loc. cit.} in that we have moved the Cartan elements
from the middle to the left.}}
\begin{equation}\label{ieq:pbw-Lsl}
e_0^{(m)}f_1^{(n)} = \sum_{p=0}^{\min(m,n)}
\sum_{j=0}^p (-1)^j \lp (\LM(h_1)-\partial_-)^{(p-j)}\cdot
f_1^{(n-p)} \rp
\lp \partial_+^{(j)}\cdot e_0^{(m-j)}\rp
\end{equation}
where,
\begin{itemize}
\item $\LM(X)$ denotes the operator of left multiplication by $X$.
\item $\partial_{\pm}$ are the derivations of $U(L\Lsl_2)$ given by
\[
\begin{aligned}
\partial_-(f_k)&=(k-2)f_{k+1}, & \partial_+(f_k)&=(k-1)f_{k+1}\\
\partial_-(h_k)&=kh_k, & \partial_+(h_k)&=kh_k\\
\partial_-(e_k)&=(k+2)e_{k+1}, & \partial_+(e_k)&=(k+1)e_{k+1}
\end{aligned}
\]
\end{itemize}

Using Garland's formula, one can show the following relation in $\power{U(L\Lsl_2)}{z^{-1}}$

\begin{equation}\label{ieq:series-Lsl}
\sum_{n=0}^{\infty} (-1)^n e_0^{(n)}f_1^{(n)}z^{-n} = 
P(z) \sum_{n=0}^{\infty} (\exp(\partial_-z^{-1})\cdot f_1)^{(n)}
(\exp(\partial_+ z^{-1})\cdot e_0)^{(n)}
\end{equation}
where $\ds P(z) = \exp\lp - \sum_{n=1}^{\infty} \frac{h_n}{n} z^{-n}\rp$.

We will
not prove \eqref{ieq:series-Lsl} here, since it is not needed, and
can be easily obtained from its $q$--analogue proved
in Section \ref{ssec:straight}.

\begin{rem}
Modulo the left ideal generated by $\{e_k\}_{k\in\Z}$,
the identity \eqref{ieq:series-Lsl} is due to Chari 
\cite[Prop.~1.1]{chari-integrable}. This is the truncation
on highest weight space result mentioned in Section \ref{issec:q-diff}.
\end{rem}

One can deduce the rationality of $P(z)$ from \eqref{ieq:series-Lsl} as follows. 
It is clear that on a \fd representation of $L\Lsl_2$, both
infinite sums in \eqref{ieq:series-Lsl} are finite. Further, noting that 
\[
\exp(\partial_-z^{-1})\cdot f_1 = f_1 - f_2z^{-1},
\text{ and }\exp(\partial_+z^{-1})\cdot e_0 = e(z) = \sum_{r=0}^{\infty}
e_rz^{-r}
\]
the rationality of these series on finite--dimensional
representations implies that of $P(z)$.

\subsection{} \label{ss:sketch end}
Our argument for the rationality of the series $\PP^\pm(z)$ defined by \eqref{ieq:diff} 
is essentially this, except in the $q$--setting.
Namely, in Proposition \ref{pr:straight}, we show that the following
equation holds in $\power{\qloopsl{2}}{z^{-1}}$

\begin{equation}\label{ieq:series-qLsl}
\begin{aligned}
\sum_{n=0}^{\infty} (-1)^n q^{n^2}E_0^{(n)}F_1^{(n)}K^{-n}z^{-n} =\hspace*{2.5in} \\
\hspace*{0.5in} 
\sum_{\ell=0}^{\infty} (-1)^{\ell} q^{\ell^2}K^{-\ell} z^{-\ell}\cdot
\PP^+(q^{-2\ell}z)\lp F_1-q^{2\ell+2}F_2z^{-1}\rp^{(\ell)}
E^+(q^{-2\ell}z)^{(\ell)}
\end{aligned}
\end{equation}
where $\{E_k,F_k\}_{k\in\Z}$ are the loop generators of $\qloopsl{2}$,
$x^{(n)} = \frac{x^n}{[n]!}$, and $E^+(z)=\sum_{r=0}^{\infty} E_rz^{-r}$. 

The proof of \eqref{ieq:series-qLsl} relies on the analogue of \eqref{ieq:pbw-Lsl}
for $\qloopsl{2}$ obtained by Chari--Pressley \cite[Lemma 5.1]{cp-rootofunity}.
Parts \eqref{iit:rat} and \eqref{iit:lim} of Theorem \ref{ithm:main}
are direct consequences
of this, combined with the rationality of half--currents
on finite--dimensional representations \cite{beck-kac,hernandez-drinfeld-coproduct,sachin-valerio-2}.
Part \eqref{iit:lattice} is proved in Section \ref{ssec:pf-3}
by (a) checking that both sides
have the same commutation relations with raising/lowering operators
of $\qloopsl{2}$ in Section \ref{ssec:pf-3-comm}, and (b) 
a direct verification of it on the subspace $\Ker(E_{-1})$
in Sections \ref{ssec:pf-3I} and \ref{ssec:pf-3II}.

\subsection{Outline of the paper} 
This paper is organized as follows. In Section
\ref{sec: notations} we fix some notations for affine Lie algebras and the corresponding
extended affine Weyl and braid groups. Section \ref{sec: qaffine} contains the
definitions of $\qloopkm$ and $\qloop$ and reviews Beck's isomorphism.
In Section \ref{sec:main}, we state the main theorem(Theorem \ref{thm:cp-series})
and prove some of its corollaries. Theorem \ref{thm:cp-series} is proved in Sections
\ref{sec:proof} and \ref{ssec:pf-3}. We review the definition of Nakajima quiver varieties,
and compute the action of the lattice operators on their equivariant $K$--theory in Section
\ref{sec:NQV}.

\subsection{Acknowledgements}
SG was supported through the Simons foundation
collaboration grant 526947, and VTL by the NSF grant DMS--2302568.

\section{Background and notations}\label{sec: notations}

In this section we set up some notations for the (untwisted) affine Lie 
algebra $\wh{\g}$ corresponding to a simple Lie algebra $\g$. We
follow \cite{kac} for this section.

\subsection{Simple Lie algebras}\label{ssec: simple-lie-algebras}

Let $\bfA=(a_{ij})_{i,j\in\bfI}$ be a Cartan matrix of finite type,
and $D=\mathrm{diagonal}(d_i)_{i\in\bfI}$ symmetrizing integers. Thus,
$a_{ii}=2$; $a_{ij}\leq 0$ for $i\neq j$; and $D\bfA$ is a symmetric,
positive--definite matrix. We will assume that $\bfA$ is indecomposable,
and $\gcd(d_i)_{i\in\bfI}=1$.

Let $(\h,\{\alpha_i\}_{i\in\bfI},\{h_i\}_{i\in\bfI})$ denote a realization of $\bfA$. 
Thus, $\h$ is $|\bfI|$--dimensional vector space over $\C$, $\{h_i\}\subset\h$
is a basis of $\h$ and $\{\alpha_i\}\subset\h^*$ the basis uniquely
determined by $\alpha_i(h_j)=a_{ji}$. Let $(\cdot,\cdot)$ denote
the symmetric, bilinear, positive--definite form on $\h^*$ (resp. $\h$)
whose matrix in the basis $\{\alpha_i\}$ (resp. $\{h_i\}$)
is $D\bfA$ (resp. $\bfA D^{-1}$). Let $\nu:\h^*\to\h$ denote
the isomorphism resulting from $(\cdot,\cdot)$. Thus, $\nu(\alpha_i)=d_ih_i$.
For $\gamma\in\h^*\setminus\{0\}$, denote $\gamma^{\vee} = 2\nu(\gamma)/(\gamma,\gamma)$
and let $s_{\gamma}\in\GL(\h^*)$ 
(resp. $\GL(\h)$) be the reflection defined
by $s_{\gamma}(\beta) = \beta-\beta(\gamma^{\vee})\gamma,\ \forall\ 
\beta\in\h^*$ (resp.
$s_{\gamma}(h) = h - \gamma(h)\gamma^{\vee},\ \forall\ h\in\h$).

Let $W\subset\GL(\h^*)$ (resp. $\GL(\h)$) denote the Weyl group, that is,
the group generated by simple reflections $\{s_i=s_{\alpha_i}\}_{i\in\bfI}$.
It is well known that $W$ is a Coxeter group, that is,
it admits the following presentation, where $m_{ij}=2,3,4,$ or $6$, if
$a_{ij}a_{ji}=0,1,2,$ or $3$ respectively.

\[
s_i^2=1,\ \forall\ i\in\bfI, \quad\text{and}\quad
\underbrace{s_is_js_i\ldots}_{m_{ij}\text{ terms}} = 
\underbrace{s_js_is_j\ldots}_{m_{ij}\text{ terms}},
\ \forall\ i\neq j\in\bfI
\]

Let $R\subset\h^*$ denote the set of roots, $R_+\subset R$
the set of positive roots and $\theta\in R_+$ the unique longest root.
Let $Q = \Z R\subset \h^*$ denote the root lattice, and
$P := \{\gamma\in\h^* : \gamma(h_i)\in\Z,\ \forall\ i\in\bfI\}$
the weight lattice. Similarly, $Q^{\vee}\subset P^{\vee}\subset\h$
denote the coroot and coweight lattices respectively.
For each $i\in\bfI$, let $\varpi_i^{\vee}\in P^{\vee}$ denote
the $i^{\text{th}}$ fundamental coweight, given by $\alpha_j(\varpi_i^{\vee})=\delta_{ij}$.

Let $\g$ denote the finite--dimensional, simple Lie algebra
over $\C$ associated to the Cartan matrix $\bfA$. 
Thus, $\g$ admits a presentation on the set of generators
$\{e_i,f_i,h_i\}_{i\in\bfI}$ with the following relations:
\begin{enumerate}
\item $[h_i,h_j]=0$
\item $[h_i,e_j] = a_{ij}e_j \aand [h_i,f_j]=-a_{ij}f_j$
\item $[e_i,f_j]=\delta_{ij} h_i$
\item For $i\neq j$,
\[
\ad(e_i)^{1-a_{ij}}(e_j) = 0 \aand \ad(f_i)^{1-a_{ij}}(f_j)=0
\]
\end{enumerate}

We continue to denote by $(\cdot,\cdot)$ the non--degenerate,
symmetric, bilinear and invariant form on $\g$, extending the
same on $\h$ by setting $(e_i,f_j)=\delta_{ij} d_i^{-1}$.

\subsection{Affine Lie algebras}\label{ssec: affine-lie-algebras}

Let $L\g := \g[z,z^{-1}]$ be the loop algebra
of $\g$. The untwisted affine Lie algebra $\wh{\g}$ associated
to $\g$ is defined as:
\[
\wh{\g} := \lp L\g \oplus \C c\rp \rtimes \C d
\]
with the following Lie bracket:
\begin{gather*}
[x(k),y(l)]=[x,y](k+l) + mk\delta_{k,-l}(x,y)c \\
\ad(c) = 0 \aand [d,x(k)]=kx(k)
\end{gather*}
for every $x,y\in \g$, $k,l\in \Z$. Here, $m=(\theta,\theta)/2$ and 
we have used the notation $x(k) = x\otimes z^k\in L\g$.\\

Let $\wh{\h} = \h\oplus\C c\oplus \C d$ be the Cartan
subalgebra of $\wh{\g}$. We extend the inner product on $\h$
to $\wh{\h}$ by declaring $(d,\h) = (c,\h) = (d,d) =
(c,c)=0$ and $(c,d) = m^{-1}$. The corresponding isomorphism
$\wh{\nu} : \wh{\h}^*\to \wh{\h}$ is then given by:
\[
\wh{\nu}|_{\h^*}=\nu \qquad \wh{\nu}(\Lambda)=md \qquad
\wh{\nu}(\delta)=mc
\]
where $\Lambda,\delta\in \wh{\h}^*$ are linear forms dual 
to $c,d$ respectively. The affine root system $\wh{R}$
then becomes:
\[
\wh{R} = \{\alpha+n\delta : \text{ either } \alpha\in R, 
n\in \Z \text{ or } \alpha=0, n\in \Z^{\times}\}
\]
We choose the following base of $\wh{R}$:
\[
\wh{\Delta} := \{\alpha_i : i\in \bfI\} \cup \{\alpha_0
:= \delta-\theta\}
\]

\subsection{Gabber-Kac isomorphism}
Let $\wh{\bfI} := \bfI\cup\{0\}$ and let $\wh{\bfA}=(a_{ij})_{i,j\in\wh{\bfI}}$
be the affine Cartan matrix, where
\[
a_{0j} = -\alpha_j(\ex{\theta}) \qquad
a_{i0} = -\theta(\ex{\alpha_i}) \qquad
a_{00} = 2
\]

Let $\g(\wh{\bfA})$ be the Kac--Moody Lie algebra associated
to $\wh{\bfA}$. Then, $\wh{\g}$ and $\g(\wh{\bfA})$ are isomorphic,
with the isomorphism given as follows. 

For a realization of $\wh{\bfA}$ we take $\wh{\h}$ as in the previous section,
$\{h_i\}_{i\in\wh{\bfI}}\subset\wh{\h}$,
where $h_0=c-\theta^{\vee}$, and $\{\alpha_i\}_{i\in\wh{\bfI}}$ where
$\alpha_0=\delta-\theta$. Thus, the Cartan subalgebras of $\wh{\g}$
and $\g(\wh{\bfA})$ are assumed to be identified.

Let $e_\theta\in\g_\theta$ and $f_\theta\in\g_{-\theta}$ be chosen
so that $(e_\theta,f_\theta) = 1/m$. Let $\{\mathfrak{e}_i,\mathfrak{f}_i\}
_{i\in\wh{\bfI}}$ denote the Chevalley--type generators of
$\g(\wh{\bfA})$.
For any $c\in\nC$, the following assignment $\GK_c : \g(\wh{\bfA})\to\wh{\g}$ 
extends to an isomorphism
of Lie algebras: $\GK_c(\mathfrak{e}_i) = e_i(0)$, $\GK_c(\mathfrak{f}_i)
=f_i(0)$ for $i\in\bfI$,
\[
\GK_c(\mathfrak{e}_0) = cf_{\theta}(1) \qquad\text{and}\qquad 
\GK_c(\mathfrak{f}_0) = c^{-1}e_{\theta}(-1)\ .
\]

\subsection{Affine Weyl group}\label{ssec: affine-weyl-group}
 
Let $W_{\aff} \subset 
GL(\wh{\h}^*)$ be defined as the subgroup generated by 
the reflections $\{s_i : i\in\wh{\bfI}\}$. 
This is again a Coxeter group, and we have an isomorphism
$W_{\aff}\isom W\ltimes Q^{\vee}$ obtained as follows.

Define $t_{\ex{\theta}}
\in W_{\aff}$ by $t_{\ex{\theta}} = s_0s_{\theta}$, and
for any $w\in W$, set $t_{w(\ex{\theta})} := wt_{\ex{\theta}} w^{-1}$.
One then shows that $\{t_{w(\ex{\theta})}\}$ commute, and
generate an abelian, normal subgroup $W'\subset W_{\aff}$, which
is isomorphic to $Q^{\vee}$.
Moreover, $W_{\aff}$ is generated by $W$ and $W'$, and
$W\cap W'=\{1\}$, thus establishing the isomorphism
$W_{\aff} \cong W \ltimes \ex{Q}$.

\subsection{Extended affine Weyl group}\label{ssec: extended-weyl-group}
Define the extended
affine Weyl group as
$\ds W_{\aff}^e := W\ltimes \ex{P}$. This is a larger group of symmetries
of the affine root system $\wh{R}$, where the $W$ action on $\wh{\h}^*$ is
as before, and for $x\in \ex{P}$ and $\xi\in\wh{\h}^*$, we have
\[
t_x(\xi) = \xi + m\xi(c) \nu^{-1}(x) - \lp \xi(x) + m\xi(c)\frac{(x,x)}{2}\rp\delta
\]
It is an easy exercise to check that the above $W_{\aff}^e$--action
preserves the set of affine roots $\wh{R}$.
For $w\in W_{\aff}^e$ define $l(w)$ to be the number of 
positive (affine) roots mapped to negatives by $w$. Let $\Pi$
be the set of elements of $W_{\aff}^e$ of length $0$. Clearly
the elements of $\Pi$ act as automorphisms of the (affine) Dynkin
diagram of $\wh{\g}$:
\[
\tau(\alpha_i) = \alpha_{\tau{i}},\ \forall\ i\in \wh{\bfI}
\]

The following theorem is well known (see e.g, 
\cite[Chapter 2, \S 6]{bourbaki-lie456}).

\begin{thm}\label{thm: extended-weyl-group}
$\Pi$ is a subgroup of $W_{\aff}^e$ isomorphic to $\ex{P}/
\ex{Q}$. There is a bijection between $\Pi\setminus \{1\}$
and the set of minuscule coweights:
\[
\Pi\setminus \{1\} \leftrightarrow \bfJ := 
\{i\in \bfI :  \theta(\ex{\varpi_i})=1\}
\]
Moreover we have an isomorphism $W_{\aff}^e \cong \Pi 
\ltimes W_{\aff}$.
\end{thm}


\subsection{Braid groups}\label{ssec: braid-gp}

Using the length function on $W_{\aff}^e$ one can define the 
associated braid group $\ebaff$ as:
\begin{equation}\label{eq: braid-gp}
\ebaff := \la\left. T_w : w\in W_{\aff}^e \right|
T_uT_v = T_{uv} \text{ if } l(uv) = l(u) + l(v)\ra
\end{equation}

Corresponding to the two incarnations $W_{\aff}^e = W \ltimes 
P^{\vee} = \Pi\ltimes W_{\aff}$ we have the following two 
presentations of $\ebaff$ (see \cite[\S 3.3]{macdonald-affine}).
Here, $m_{ij}=2,3,4,6,$ or $\infty$, if $a_{ij}a_{ji}
=0,1,2,3,$ or $4$.

\begin{prop}\label{prop: braid-gp}
\hfill
\begin{enumerate}
\item $\ebaff$ is generated by $U_{\tau} (\tau\in \Pi)$ and 
$T_i (i\in \wh{\bfI})$ subject to the following relations:
\begin{gather}
U_{\tau}U_{\tau'} = U_{\tau+\tau'}\\
\braid{T}{i}{j}{m_{ij}} = \braid{T}{j}{i}{m_{ij}}\\
U_{\tau}T_iU_{\tau}^{-1} = T_{\tau(i)}
\end{gather}

\item Let $Y_i := T_{t_{\ex{\varpi_i}}} \in \ebaff$. Then $\ebaff$ is 
generated by $\{T_i, Y_i\}_{i\in \bfI}$ subject to the following relations:
\begin{gather}
\braid{T}{i}{j}{m_{ij}} = \braid{T}{j}{i}{m_{ij}} \label{eq: TT}\\
Y_iY_j = Y_jY_i \label{eq: YY}\\
T_i^{-1}Y_iT_i^{-1} = Y_i^{-1} \prod_{j\neq i} Y_j^{-a_{ij}} 
\label{eq: TY1}\\
T_iY_j = Y_jT_i \text{ for } j\neq i \label{eq: TY2}
\end{gather}
We will use the notation $L_i := T_{t_{\ex{\alpha_i}}} \in \baff \subset \ebaff$,
where $\baff$ denotes the affine braid group. Note that $\baff$
is a normal subgroup of $\ebaff$.
\end{enumerate}
\end{prop}

\begin{example}\label{ex:sl2} 
For $\g=\Lsl_2$, the affine braid group $\baff$ is the
free group on two generators $\{T_0,T_1\}$, and $L=T_0T_1$. The diagram
automorphism $\sigma$ swaps $0$ and $1$.
The extended affine braid group $\ebaff$ has the
following two presentations:
\[
\ebaff=\langle U_{\sigma},T_0,T_1 | U_{\sigma}^2=1,\ 
U_\sigma T_0 U_\sigma^{-1} = T_1 \rangle
\isom \langle T, Y | T^{-1}YT^{-1} = Y^{-1}\rangle
\]
where, the isomorphism maps $U_{\sigma}\mapsto YT^{-1}$, $T_1\mapsto T$, and hence
$T_0 = \tau T_1\tau \mapsto Y^2T^{-1}$. Thus, the coroot lattice element
$L\mapsto Y^2$, $Y$ being the coweight lattice element.
\end{example}

\noindent The following will be crucial in carrying out a rank one reduction argument.\\

\begin{cor}\label{cor: rank1-formula}
For each $i\in\bfI$, we have:
\[
L_i = \lp\Ad(Y_i)\cdot T_i \rp  T_i
\]
\end{cor}

\begin{pf}
Note that, upon taking inverses, and using $L_i = \prod_j Y_j^{a_{ij}}$,
the relation \eqref{eq: TY1} becomes
\[
T_iY_i^{-1}T_i = Y_i^{-1}L_i
\]
\end{pf}

\section{Quantum loop algebra}\label{sec: qaffine}

In this section we review the two presentations of the quantum loop algebra
associated to $\g$, the action of the extended affine braid group, and Beck's
isomorphisms.

\subsection{Kac--Moody presentation}
\label{ssec: km-presentation}

Let $q\in\C^\times$ be of infinite order and set $q_i=q^{d_i}$ for any $i\in\bfI$.
We use the standard notation for Gaussian
integers
\begin{gather*}
[n]_i = \frac{q_i^n - q_i^{-n}}{q_i-q_i^{-1}}\\[.5ex]
[n]_i! = [n]_i[n-1]_i\cdots [1]_i\qquad
\qbin{n}{k}{i} = \frac{[n]_i!}{[k]_i![n-k]_i!}
\end{gather*}

The quantum loop algebra $\qloopkm$ is the $\C$--algebra generated by elements
$\{\K_i^{\pm 1},\E_j, \F_j\}_{i\in \bfI, j\in \wh{\bfI}}$ subject to the relations
\begin{itemize}
\item[(KM1)] $[\K_i,\K_j]=0$ for every $i,j\in \bfI$.
\item[(KM2)] For every $i\in \bfI$ and $j\in \wh{\bfI}$
\[
\K_i \E_j \K_i^{-1} = q_i^{a_{ij}} \E_j\aand
\K_i \F_j \K_i^{-1} = q_i^{-a_{ij}} \F_j\]
\item[(KM3)] For each $i,j\in \wh{\bfI}$ we have
\[
[\E_i,\F_j] = \delta_{ij} \frac{\K_i-\K_i^{-1}}{q_i-q_i^{-1}}
\]
where $\K_0 = \K_{\theta}^{-1}$ and $d_0=m$. Here, if $\theta=\sum_{i\in\bfI} n_i\alpha_i$,
then $\K_{\theta} := \prod_{i\in\bfI} \K_i^{n_i}$.
\item[(KM4)] For $i\neq j \in \wh{\bfI}$, we have
\begin{gather*}
\sum_{s=0}^{1-a_{ij}} (-1)^s \qbin{1-a_{ij}}{s}{q_i} 
\E_i^{1-a_{ij}-s}\E_j\E_i^s = 0 \\
\sum_{s=0}^{1-a_{ij}} (-1)^s \qbin{1-a_{ij}}{s}{q_i} 
\F_i^{1-a_{ij}-s}\F_j\F_i^s = 0 
\end{gather*}
\end{itemize}

\subsection{Quantum Weyl group}
\label{ssec: q-weyl-group}

Let $V$ be a type I, integrable representation of $\qloopkm$.
That is, the eigenvalues of $\K_i$ lie in $q_i^{\Z}$
and $\E_i,\F_i$ $(i\in\wh{\bfI})$ act locally
nilpotently
\footnote{meaning, for every $v\in V$, there exists $N\gg 0$,
such that $\E_i^Nv = \F_i^Nv = 0$.
This holds, for instance, when $V$ is finite--dimensional.}.
Let $\HH_i$ denote the unique semisimple operator on $V$ with $\Z$ eigenvalues
such that $\K_i=q_i^{\HH_i}$. In this paper we solely work with
type I representations, which will be assumed throughout.

Then, the affine braid group $\baff$ acts naturally on $V$
\cite{kirillov-reshetikhin, lusztig-canonicalII, lusztig-book, saito-PBW, soibelman}. In more detail,
we have a group homomorphism $\qW_V : \baff \to \GL(V)$, where
$\qW_V(T_i)$ is given by the following triple $q$--exponential
formula:
\begin{equation}\label{eq: triple-ex}
\qS_i := \exp_{q_i^{-1}}(q_i^{-1}\E_i\K_i^{-1})
\exp_{q_i^{-1}}(-\F_i)\exp_{q_i^{-1}}(q_i\E_i\K_i)
q_i^{\HH_i(\HH_i+1)/2}
\end{equation}
Thus, $\qW_V(T_i) = \pi_V(\qS_i)$, where $\pi_V:\qloopkm\to\End(V)$
is the action homomorphism. Note that the integrability hypothesis
on $V$ is necessary to make sense of $\pi_V(\qS_i)$.
A complete proof of this result can be found in \cite[Ch. 39]{lusztig-book}.

\subsection{Braid group action on $\qloopkm$}\label{ssec: auto}
The braid group action $\qW_V$ on every \fd representation
$V\in\Rloopkm$ lifts to an action on $\qloopkm$
\cite{lusztig-qdef,lusztig-book}. More precisely,
there is $\qW:\baff\to\Aut(\qloopkm)$ (see, \eg \cite[Ch. 37]{lusztig-book})
such that
\begin{equation}\label{eq:qw-consistent}
	\qW_V(g) (a\cdot (\qW_V(g^{-1})(v))) = \qW(g)(a) \cdot v
\end{equation}
for every $V\in\Rloopkm,\ v\in V,\ a\in\qloopkm$ and $g\in\baff$.

Thus, for the generator $T_i$ of $\baff$, $\qW(T_i) = \Ad(\qS_i)$.
One needs to verify that the latter makes sense, since $\qS_i$ are
elements of a completion of $\qloopkm$. This follows from
the explicit formulae for the action of $T_i$ on 
the generators $\{\E_j, \F_j, \K_j\}_{j\in \wh{\bfI}}$ (see \cite[Ch. 36]{
lusztig-book}).

\begin{align}\label{eq: braid-KM1}
T_i(\E_i) &= -\F_i\K_i
&T_i(\F_i) &= -\K_i^{-1}\E_i
&T_i(\K_j) &= \K_j\K_i^{-a_{ij}}
\end{align}
and for $i\neq j \in \wh{\bfI}$
\begin{gather}
T_i(\E_j) = \sum_{s=0}^{-a_{ij}} (-1)^{s-a_{ij}} q_i^{-s} 
\E_i^{(-a_{ij}-s)}\E_j\E_i^{(s)} \label{eq: braid-KM2-1}\\
T_i(\F_j) = \sum_{s=0}^{-a_{ij}} (-1)^{s-a_{ij}} q_i^{s} 
\F_i^{(s)}\F_j\F_i^{(-a_{ij}-s)} \label{eq: braid-KM2-2}
\end{gather}

We extend this action of $\baff$ on $\qloopkm$ to an action of 
$\ebaff$ by defining (see Proposition \ref{prop: braid-gp})
\begin{align}
U_{\tau}\E_i &= \E_{\tau(i)}
&U_{\tau}\F_i &= \F_{\tau(i)}
&U_{\tau}\K_i &= \K_{\tau(i)}
\end{align}
By a little abuse of notation we continue to denote it by
$\qW:\ebaff\to\Aut(\qloopkm)$.

\begin{lem}\label{lem:conj-id}
Let $V$ be a finite--dimensional representation of $\qloopkm$.
Then, for every $g\in\ebaff$ and $x\in\baff$, the following
equation holds
\begin{equation}\label{eq:conj-id}
\qW_V(gxg^{-1}) = \qW_{g^*V}(x)
\end{equation}
where $g^*V = \qW(g)^*(V)$ is the pull--back representation under
the algebra isomorphism $\qW(g)$.
\end{lem}

\begin{pf}
Let $\pi_V : \qloopkm\to\End(V)$ denote the action homomorphism.
When both $g,x\in\baff$, the equation \eqref{eq:conj-id} is
a consequence of \eqref{eq:qw-consistent}, which can be seen as follows.

\[
\qW_V(g)\circ\pi_V(a)\circ\qW_V(g^{-1}) = \pi_V(\qW(g)(a))
= \pi_{g^*V}(a)
\]

Now take $a\in\qloopkm$ so that $\pi_V(a) = \qW_V(x)$. That such an $a$
exists (depending on $V$) follows from the fact that $\qW_V(T_i)$
is defined as $\pi_V(\qS_i)$ given in \eqref{eq: triple-ex}.

Thus, it suffices to consider the case when $g=U_{\sigma}$ and $x=T_i$,
for some $\sigma\in\Pi$, $i\in\bfI$. By Proposition \ref{prop: braid-gp},
$U_{\sigma} T_i U_{\sigma}^{-1} = T_{\sigma(i)}$ and \eqref{eq:conj-id} becomes
$\qW_V(T_{\sigma(i)}) = \qW_{U_{\sigma}^*V}(T_i)$. This
follows from the expression \eqref{eq: triple-ex} of
$\qS_i$ and the fact that $U_\sigma(X_i) = X_{\sigma(i)}$, for $X=\E,\F$
or $\HH$. The lemma is proved.
\end{pf}

\subsection{Rank $1$ reduction}\label{ssec:rk1-red}
Let $i\in\bfI$ and let $Y_i' := Y_iT_i^{-1}$. By \cite[Prop. 3.8]{beck-braid}, the following is an
injective algebra homomorphism $\vembed_i : \qlsl{i} \to \qloopkm$.

\begin{equation}
\begin{aligned}
\vembed_i(\E_1) &= \E_i
&\vembed_i(\F_1) &= \F_i \\
\vembed_i(\E_0) &= Y_i'(\E_i)
&\vembed_i(\F_0) &= Y_i'(\F_i)
\end{aligned}
\end{equation}

\begin{prop}\label{pr:rk1-red}
For each $i\in\bfI$ and $V\in\Rloopkm$ we have
\[
\qW_V(L_i) = \qW_{\vembed_i^*(V)}(L)
\]
where $L = T_0T_1$ is the lattice operator for $\Lsl_2$.
\end{prop}

\begin{pf}
Let us write $\Phi(\E,\F,\H)$ for the triple exponential:
\[
\Phi(\E,\F,\H) = 
\exp_{q^{-1}}(q^{-1}\E\K^{-1})
\exp_{q^{-1}}(-\F)\exp_{q^{-1}}(q\E\K)
q^{\HH(\HH+1)/2}
\]
By definition of $\vembed_i$, $\qW_V(T_i) = \qW_{\vembed_i^*(V)}(T_1)$.
Thus, by
Corollary \ref{cor: rank1-formula} and $\Ad(Y_i')(T_i) = \Ad(Y_i)(T_i)$,
it suffices to verify that
\[
\qW_V(\Ad(Y_i')(T_i)) = \qW_{\vembed_i^*(V)}(T_0)
\]

Using equation \eqref{eq:conj-id} we get that

\begin{align*}
\qW_V(\Ad(Y_i')(T_i)) &= \Phi(\E_i,\F_i,\HH_i) \text{ on } (Y_i')^*(V) \\
&= \Phi(Y_i'(\E_i),Y_i'(\F_i),Y_i'(\HH_i)) 
= \Phi(\vembed_i(\E_0),\vembed_i(\F_0),\vembed_i(\HH_0)) \\
&= \qW_{\vembed_i^*(V)}(T_0)
\end{align*}
as claimed.
\end{pf}

\subsection{Loop presentation \cite{drinfeld-yangian-qaffine}}
\label{ssec: loop-presentation}
Let $\qloop$ be the unital, associative algebra over $\C$, 
generated by $\{E_{i,k}, F_{i,k}, \cg^{\pm}_{i,\pm\ell}\}_
{i\in \bfI, k\in \Z, \ell\in\Z_{\geq 0}}$
subject to the following relations:
\begin{itemize}
\item[(QL1)] For $i,j\in\bfI$ and $r,s\in\Z_{\geq 0}$
\[[\cg^{\pm}_{i,\pm r},\cg^{\pm}_{j,\pm s}]=0 = 
[\cg^{\pm}_{i,\pm r},\cg^{\mp}_{j,\mp s}]
\]
Moreover, $\cg^+_{i,0}\cg^-_{i,0} = 1$ for each $i\in\bfI$.
\item[(QL2)] For any $i,j\in\bfI$ and $k\in \Z$,
\[
\Ad(\cg^{\pm}_{i,0})\cdot E_{j,k} = q_i^{\pm a_{ij}}E_{j,k}\qquad
\Ad(\cg^{\pm}_{i,0})\cdot F_{j,k} = q_i^{\mp a_{ij}}F_{j,k}
\]

\Omit{\item[(QL3)] For any $i,j\in\bfI$ and $r\in\Z^\times$,
\[[H_{i,r}, E_{j,k}]= \frac{[ra_{ij}]_{q_i}}{r} E_{j,r+k}
\qquad
[H_{i,r}, F_{j,k}]= -\frac{[ra_{ij}]_{q_i}}{r}F_{j,r+k}\]
}

\item[(QL3)] For any $i,j\in\bfI$, $\eta\in\{\pm\}$ and $\ell\in\Z$,
\begin{align*}
\cg^\eta_{i,k+1}E_{j,l} - q_i^{a_{ij}}E_{j,l}\cg^{\eta}_{i,k+1} 
&= q_i^{a_{ij}}\cg^{\eta}_{i,k}E_{j,l+1} - E_{j,l+1}\cg^{\eta}_{i,k}\\
\cg^{\eta}_{i,k+1}F_{j,l} - q_i^{-a_{ij}}F_{j,l}\cg^{\eta}_{i,k+1} 
&= q_i^{-a_{ij}}\cg^{\eta}_{i,k}F_{j,l+1} - F_{j,l+1}\cg^{\eta}_{i,k}
\end{align*}
for each $k\in\Z_{\geq 0}$ if $\eta=+$, and $k\in\Z_{<0}$ if
$\eta=-$. Here, and in the future, we follow the convention
that $\cg^{\pm}_{i,\mp r}=0$ for each $r\in\Z_{\geq 1}$.
\item[(QL4)] For $i,j\in\bfI$ and $k,l\in \Z$
\begin{align*}
E_{i,k+1}E_{j,l} - q_i^{a_{ij}}E_{j,l}E_{i,k+1} &= q_i^{a_{ij}}E_{i,k}E_{j,l+1} - E_{j,l+1}E_{i,k}\\
F_{i,k+1}F_{j,l} - q_i^{-a_{ij}}F_{j,l}F_{i,k+1} &= q_i^{-a_{ij}}F_{i,k}F_{j,l+1} - F_{j,l+1}F_{i,k}
\end{align*}
\item[(QL5)] For $i,j\in\bfI$ and $k,l\in \Z$
\[[E_{i,k}, F_{j,l}] = \delta_{ij} \frac{\cg^+_{i,k+l} - \cg^-_{i,k+l}}
{q_i - q_i^{-1}}\]
\item[(QL6)] Let $i\not= j\in\bfI$ and set $m = 1-a_{ij}$. For every $k_1,\ldots, k_m\in \Z$
and $l\in \Z$
\begin{gather*}
\sum_{\pi\in \Sym_m} \sum_{s=0}^m (-1)^s\bin{m}{s}{q_i}
E_{i,k_{\pi(1)}}\cdots E_{i,k_{\pi(s)}} E_{j,l}E_{i,k_{\pi(s+1)}}\cdots E_{i,k_{\pi(m)}} = 0\\
\sum_{\pi\in \Sym_m} \sum_{s=0}^m (-1)^s\bin{m}{s}{q_i}
F_{i,k_{\pi(1)}}\cdots F_{i,k_{\pi(s)}} F_{j,l}F_{i,k_{\pi(s+1)}}\cdots F_{i,k_{\pi(m)}} = 0
\end{gather*}
\end{itemize}

Later we will need another system of generators of the maximal
commutative subalgebra $U^0_q(L\g)$ of $\qloop$, denoted by 
$\{\K_i,H_{i,k}\}_{i\in\bfI, k\in\Z_{\neq 0}}$, defined by
the following equation:

\begin{equation}\label{eq: psi-h}
\cg_i^{\pm}(z):= \sum_{k\geq 0} \cg^{\pm}_{i,\pm k} z^{\mp k}
 = \K_i^{\pm 1}\exp\lp\pm(q_i-q_i^{-1})\sum_{r\geq 1} H_{i,\pm r}
 z^{\mp r}\rp
\end{equation}

Let $E_i(w)=\sum_{\ell\in\Z} E_{i,\ell} z^{-\ell}$ and $F_i(w)=\sum_{\ell\in\Z}
F_{i,\ell}z^{-\ell}$. Relation (QL3) is equivalent to the following
identity in $\qloop[\![z^{\mp 1}, w,w^{-1}]\!]$

\begin{equation}\label{eq:psi-EF}
\begin{split}
\Ad(\cg_i^{\pm}(z))\cdot E(w) &= \frac{q_i^{a_{ij}}z-w}{z-q^{a_{ij}}w} E(w) \\
\Ad(\cg_i^{\pm}(z))^{-1} \cdot F(w) &= \frac{q_i^{a_{ij}}z-w}{z-q^{a_{ij}}w} F(w)
\end{split}
\end{equation}

\subsection{Shift automorphisms}\label{ssec:shift}
$\qloop$ admits a $1$--parameter group of algebra automorphisms,
denoted by $\tau_\zeta\in\Aut(\qloop)$, $\zeta\in\nC$, given as:
\[
\tau_\zeta(X_k) = \zeta^k X_k,\ \text{ where } X \text{ is one of }
E_i,F_i,\cg^{\pm}_i\ \ (i\in\bfI).
\]

In terms of the formal series defined above, we have
$\tau_\zeta(X(z)) = X(\zeta^{-1}z)$, for $X=E_i,F_i,\cg_i^{\pm}$.
For a representation $V$ of $\qloop$, we denote the pull--back
representation $\tau_{\zeta}^*(V)$ by $V(\zeta)$.

\subsection{Beck isomorphism}\label{ssec:beck}
We now describe Beck's isomorphisms between $\qloopkm$ and $\qloop$.
These depend on a choice of a sign $\bsign :\bfI\to\{\pm 1\}$
such that, for every $i\neq j$, $a_{ij}\neq 0 \Rightarrow \bsign(i)\bsign(j)=-1$.
As our Dynkin diagram is connected, there are two such choices.\\

Following \cite[Thm. 4.7]{beck-braid}, let us
define $\biso{\bsign}:\qloop \to \qloopkm$, as follows.
For each $i\in\bfI$ and $r\in\Z$, we have: 
$\biso{\bsign}(\cg^{\pm}_{i,0}) := \K_i^{\pm 1}$, and

\begin{equation}\label{eq: higher-modes}
\biso{\bsign}(E_{i,r}) := \bsign(i)^rY_i^{-r}\E_{i} \aand 
\biso{\bsign}(F_{i,r}) := \bsign(i)^rY_i^r\F_{i}
\end{equation}

\begin{thm}\cite[Thm. 4.7]{beck-braid}\label{thm: loop-presentation}
For each choice of the sign $\bsign:\bfI\to\{\pm 1\}$ as above,
the assignment $\biso{\bsign}$ extends to an algebra isomorphism
between $\qloop$ and $\qloopkm$.
\end{thm}

Thus, given a finite--dimensional representation $\pi_V : \qloop\to
\End(V)$, we obtain an action $\qW_{V,\bsign}:\baff\to\GL(V)$
depending on the choice of a sign $\bsign$. Thus, $\qW_{V,\bsign}$
is nothing but $\qW_{(\biso{\bsign}^{-1})^*(V)}$, the action of
the quantum Weyl group (see Section \ref{ssec: q-weyl-group}) 
on the representation $\pi_V\circ \biso{\bsign}^{-1}
:\qloopkm\to\End(V)$. When the representation $V$ is clear from
the context, we simply denote this action by $\lambdao$, for
the ease of notations.

\begin{rem}
Note that the two isomorphisms are related by
the shift automorphism $\tau_{-1}$. That is, we have:
\begin{equation}\label{eq:beck-shift}
\biso{-\bsign}^{-1}\circ\biso{\bsign} = \tau_{-1}
\end{equation}

For $\g=\Lsl_2$, the two isomorphisms are denoted by
$\biso{\pm}$. We have the following commutative diagram:
\begin{equation}\label{rk1-diagram}
\xy
(0,20)*{\qlsl{i}}="a";
(0,0)*{U_{q_i}(L\Lsl_2)}="b";
(50,20)*{\qloopkm}="c";
(50,0)*{\qloop}="d";
{\ar_{\biso{\bsign(i)}} "b"; "a"};
{\ar_{\vembed_i} "b"; "d"};
{\ar_{\vembed_i} "a"; "c"};
{\ar_{\biso{\bsign}} "d"; "c"};
\endxy
\end{equation}

We also remark that the isomorphism used for $\Lsl_2$ in Section \ref{sec:proof}
below is $\biso{-}$. Its inverse is given by:

\begin{equation}\label{eq:beta-inverse}
\begin{aligned}
\biso{-}^{-1}(\E_1) &= E_0
&\biso{-}^{-1}(\F_1) &= F_0
&\biso{-}^{-1}(\H_1) &= H_0 \\
\biso{-}^{-1}(\E_0) &= K^{-1}F_1
&\biso{-}^{-1}(\F_0) &= E_{-1}K
&\biso{-}^{-1}(\H_0) &= -H_0
\end{aligned}
\end{equation}

To see this, we recall the computation of the action of the lattice element 
$Y=T_{\varpi^{\vee}}$,
performed under the identification above in \cite[Lemmas 9.2, 9.3]{sachin-valerio-3}.

\[
Y(H_k) = H_k,\qquad Y(E_k) = -E_{k-1} \qquad
Y(F_k) = -F_{k+1}.
\]

Comparing with the formulae for $\biso{\bsign}(E_r)$, we conclude that
the isomorphism used is in fact $\biso{-}$.

\end{rem}

\section{Lattice operators}
\label{sec:main}

\subsection{Rationality of half--currents}

Let $\cg_i^{\pm}(z), E^{\pm}_i(z),F^{\pm}_i(z)\in\qloop[\![z^{\mp 1}]\!]$, 
$i\in\bfI$, be the generating series: $\cg_i^{\pm}(z) = \sum_{r\geq 0}
\cg^{\pm}_{i,\pm r} z^{\mp r}$, and

\begin{align*}
E^+_i(z) &= \sum_{r\geq 0}E_{i,r} z^{-r} & E^-_i(z) &= -\sum_{r< 0}E_{i,r}^\pm z^{-r}\\
F^+_i(z) &= \sum_{r\geq 0}F_{i,r} z^{-r} & F^-_i(z) &=-\sum_{r< 0}F_{i,r}^\pm z^{-r}
\end{align*}

The following is well--known \cite[\S 6]{beck-kac}, \cite[Prop. 38]{hernandez-drinfeld-coproduct},
\cite[Prop. 3.6]{sachin-valerio-2}.

\begin{prop}\label{pr:rat currents}
Let $V$ be a finite--dimensional representation of $\qloop$. Then,
the evaluations of $\cg^{\pm}_i(z) ,E^{\pm}_i(z) ,F^{\pm}_i(z)$ are
the Taylor series of rational $\End(V)$--valued functions 
$\cg_i(z),E_i(z),F_i(z)$ at $z=\infty,0$.
\end{prop}

\subsection{The series $\PP_i^{\pm}(z)$}\label{ssec:cp-series}

For any $i\in\bfI$, define $\cgn_i^\pm(z)\in\power{\qloop}{z^{\mp 1}}$ by
\begin{equation}\label{eq: psi-h-2}
\cgn_i^\pm(z)=\K_i^{\mp 1}\psi_i^\pm(z)
=
\exp\left(\pm(q_i-q_i^{-1})\sum_{r\geq 1} H_{i,\pm r}
 z^{\mp r}\right)
\end{equation}
so that $\cgn_i^+(\infty)=1=\cgn_i^-(0)$.

Let $\PP_i^{\pm}(z)\in\power{\qloop}{z^{\mp 1}}$ be the unique solution of the 
$q$--difference equation
\begin{equation}\label{eq:cp-diff}
\PP_i^{\pm}(q_i^2z) = 
\cgn^\pm_i(z) \,\PP_i^{\pm}(z)
\end{equation}
such that $\PP_i^+(\infty) = 1 = \PP_i^-(0)$.
By taking logarithm of this $q$--difference equation, and using
\eqref{eq: psi-h-2}, we obtain the following explicit formula for
$\PP_i^{\pm}(z)$, which also appeared in \cite[Lemma 3.2]{cp-rootofunity}.
\begin{equation}\label{eq:cp-explicit}
\PP_i^\pm(z)
=
\exp\lp -\sum_{n=1}^{\infty} q_i^{\pm n} \frac{H_{i,\pm n}}{[n]_i}z^{\mp n}\rp
\end{equation}

The series $\PP_i^{\pm}(z)$ were introduced in \cite{cp-qaffine}, and their classical
version in \cite{chari-integrable}. In \cite[Prop. 3.5]{cp-qaffine}, Chari--Pressley show
that they truncate on highest weight vectors, thus obtaining the ``only if" part of the
classification of irreducible, \fd representations of $\qloop$ via Drinfeld polynomials.

\subsection{Rationality of $\PP_i^\pm(z)$}\label{ssec:mainthm}

Let $V$ be a \fd $\qloop$--module. By Proposition \ref{pr:rat currents}, the series
$\cgn_i^\pm(z)\in \End(V)[\![z^{\mp 1}]\!]$ are Taylor expansions of $\End(V)$--valued
rational functions at $z=\infty,0$, taking value $1$ at these points. 
If $|q|\neq 1$, a standard argument then implies that $\PP_i^\pm(z)$ are
expansions of meromorphic functions
\[P_i^+(z):\IP^1\setminus\{0\}\to \GL(V)
\aand
P_i^-(z):\IP^1\setminus\{\infty\}\to \GL(V)\]
at $z=\infty,0$ normalised to take value $1$ at these points.
That is, $P_i^\pm(z)$ are the canonical solutions of the (regular) 
$q$--difference
equations \eqref{eq:cp-diff} at $z=\infty,0$, and are given by 
the convergent products
\[
P_i^+(z) = \cgn_i^+(z)^{-1}\cgn_i^+(q_i^{2}z)^{-1}\cdots
\qquad
P_i^-(z) = \cgn_i^-(q_i^{-2}z)\cgn_i^-(q_i^{-4}z)\cdots
\]
if $|q|>1$, and by
\[
P_i^+(z) =\cgn_i^+(q_i^{-2}z)\cgn_i^+(q_i^{-4}z)\cdots
\qquad
P_i^-(z) = \cgn_i^-(z)^{-1}\cgn_i^-(q_i^{2}z)^{-1}\cdots
\]
if $|q|<1$.

The main result of this paper is Theorem \ref{thm:cp-series} below,
which significantly strengthens this result, and is valid whenever
$q$ is of infinite order ($|q|$ could be $1$). 
It shows that $P^\pm_i(z)$ are in fact rational functions, or equivalently
that the difference equations \eqref{eq:cp-diff} have no monodromy. Moreover, the normalised
limits of $P^\pm_i(z)$ at $z=0,\infty$ are given by the lattice operators 
$L_i = T_{t_{\alpha_i^{\vee}}}$.

\begin{thm}\label{thm:cp-series}
The following holds for any $i\in\bfI$.
\begin{enumerate}
\item\label{it:rationality} 
$P^\pm_i(z):\C^\times\to \End(V)$ is a rational function.
\item\label{it:limits} There is an element $\cplim_i\in\GL(V)$ 
such that
\[z^{\HH_i}P_i^+(z) = \cplim_i P_i^-(z)\]
In particular,
\[\cplim_i = \lim_{z\to 0} z^{\HH_i}P_i^+(z)
=
\lim_{z\to\infty} z^{\HH_i}(P_i^-(z))^{-1}\]
\item\label{it:lattice} 
Let $\qloopkm$ act on $V$ via Beck's isomorphism 
(see Thm. \ref{thm: loop-presentation}). Then,
\[
\lambdao(L_i) = o(i)^{\HH_i} \K_i \cplim_i^{-1}
\]
\end{enumerate}
\end{thm}
\Omit{
\begin{rem}
The rationality of the matrix coefficients of $P_i^\pm(z)$ corresponding to joint
eigenvectors of $\{\cg_{i,\pm r}^\pm\}_{r\in\N}$ can also be deduced from the
form of the corresponding eigenvalues \cite[Prop. 1]{frenkel-reshetikhin-qchar}.
Part (1) of Theorem \ref{thm:cp-series} is a stronger statement, however, since
$\{\cg_{i,\pm r}^\pm\}_{r\in\N}$ do not act semisimply on $V$ in general.
\end{rem}}

\begin{pf}
The assertions \eqref{it:rationality} and \eqref{it:limits} are rank
$1$ statements, which are proved in Section \ref{sec:proof} below.
\eqref{it:lattice} reduces to rank $1$ by Proposition
\ref{pr:rk1-red}, and is proved in Section \ref{ssec:pf-3}.
\Omit{
Theorem \ref{thm:cp-series} for $\g=\Lsl_2$ is proved in 
Section \ref{sec:proof}. The proof is based on 
Proposition \ref{pr:straight}, which hinges on the remarkable observation
of Chari--Pressley that the series  $\PP_i^{\pm}(z)$ can be obtained,
modulo the left ideal generated by $\{E_{j,k}\}_{j\in\bfI,k\in\Z}$,
as the coefficients of the PBW straightening formula for $E_0^{(n)}F_1^{(n)}$ 
\cite[Prop.~3.5]{cp-qaffine}. 
Proposition \ref{pr:straight} below is a stronger version of this, as in
it does not involve quotienting by the aforementioned ideal.
This, and the rationality property of the
half--currents are then used in \ref{ssec:pf-1} to prove the rationality of
$P^\pm(z)$.

The rationality, and equality, of $\cg_i^\pm(z)$ on \fd representations then
implies \eqref{it:limits} (see \ref{ssec:pf-2}).

}
\end{pf}

In the remainder of this section, we give a few corollaries
of this theorem.

\subsection{Shift automorphisms}
\label{ssec:cp-shift}

Recall the definition of the shift automorphism $\tau_\zeta$ $(\zeta\in\nC)$
from Section \ref{ssec:shift} above.

\begin{cor}\label{cor:cp-shift}
For each $i\in\bfI$ and $\zeta\in\nC$,
$\tau_\zeta(L_i) = \zeta^{-H_{i,0}} L_i$. More precisely,
let $V\in\Rloop$ and $\zeta\in\nC$. Then, we have:
\[
\qW_{V(\zeta),\bsign}(L_i) = \zeta^{-H_{i,0}} \lambdao(L_i)
\]
\end{cor}

\begin{pf}
By Theorem \ref{thm:cp-series} \eqref{it:limits}, we have:
\[
\begin{aligned}
\cplim_{i,V(\zeta)} &= \lim_{z\to 0} z^{H_{i,0}} \pi_{V(\zeta)} (\PP_i^+(z))
= \lim_{z\to 0} z^{H_{i,0}} \pi_V(\PP^+_i(\zeta^{-1}z)) \\
&= \zeta^{H_{i,0}}\lim_{z\to 0} (\zeta^{-1}z)^{H_{i,0}} \pi_V(\PP^+_i(\zeta^{-1}z))
= \zeta^{H_{i,0}} \cplim_i
\end{aligned}
\]
The corollary now follows from \eqref{it:lattice} of Theorem
\ref{thm:cp-series}.
\end{pf}

\subsection{Euler transform and a formula for the lattice operators}
\label{ssec:euler}
We now give a closed form expression of the lattice operators in terms of
the elements $\{H_{j,\ell}\}_{j\in\bfI,\ell\in\Z_{\neq 0}}$ (see
equation \eqref{eq: psi-h}). These formulae were obtained for $\g=\Lsl_2$
in the formal
$\hbar$--adic setting, in \cite{sachin-valerio-3}. For $i\in\bfI$
and $r\geq 1$, define
\[\wt{H}_{i,r}= H_{i,0} + \sum_{s=1}^r (-1)^s \cbin{r}{s} \frac{s}{[s]_i} H_{i,s}\]
and let
\[
\mathcal{C}(t) := \exp\lp\sum_{n=1}^{\infty} \frac{\wt{H}_{i,n}}{n} t^n\rp
\in \power{\qloop}{t}
\]

\begin{cor}\label{cor:euler}
The evaluation of $\mathcal{C}(t)$ on a finite--dimensional representation
$V$ of $\qloop$ is a rational function of $t$, regular at $t=1$, and we
have:
\[
\lambdao(L_i) = (-\bsign(i))^{\H_i} \lim_{t\to 1} \mathcal{C}(t)
\]
where $\bsign$ is the sign chosen in the isomorphism $\biso{\bsign}:\qloop\to\qloopkm$.
\end{cor}

\begin{pf}
Since the statement is for one node only, we drop the subscript $i$
for convenience. Let $\epsilon=\bsign(i)\in\{\pm 1\}$. Define:
\[
\mathcal{A}(z) := (1-qz^{-1})^{H_0} \PP^+(z)^{-1}
\]
Note that, by Theorem \ref{thm:cp-series} \eqref{it:rationality}
and \eqref{it:limits}, the evaluation of $\mathcal{A}(z)$ on
$V$ is a rational function of $z$. Moreover, 
$\lim_{z\to 0} \mathcal{A}(z) = \lim_{z\to 0} (z-q)^{H_0} z^{-H_0}\PP^+(z)^{-1}
=(-1)^{H_0}\mathcal{K}\cplim^{-1}$.

We claim that $\mathcal{A}(z)$ and $\mathcal{C}(t)$ are related
by the following change of variables\footnote{Such substitutions
appeared in Euler's work, to turn slowly convergent
(or, sometimes even divergent) series to series with faster convergence rate.
We refer an interested reader to \cite[\S 8.2, 8.3]{hardy-divergent}
for more on this topic.}.
\[
z \mapsto q\frac{t-1}{t}\qquad t = \frac{1}{1-q^{-1}z}
\]
Note that $z\to 0$ corresponds to $t\to 1$.
The corollary will follow from Theorem \ref{thm:cp-series}
\eqref{it:lattice}, once the following identity is established.

\begin{equation}\label{eq:euler-pf}
\mathcal{A}(z)|_{z=q(t-1)/t} = \exp\lp\sum_{n=1}^{\infty}
\frac{\wt{H}_n}{n} t^n\rp
\end{equation}

Using the formula given in \eqref{eq:cp-explicit}, we have
\[
\mathcal{A}(z) = \exp\lp H_0\log(1-qz^{-1}) + \sum_{n=1}^{\infty}
\frac{H_n}{[n]} q^nz^{-n}\rp
\]

\begin{align*}
\left.\log(\mathcal{A}(z))\right|_{z=q(t-1)/t} &=  
-H_0\log(1-t) + \sum_{n=1}^{\infty}
(-1)^n \frac{H_n}{[n]} t^n (1-t)^{-n}\\
&= -H_0\log(1-t) + \sum_{n=1}^{\infty}
(-1)^n \frac{H_n}{[n]} t^n \lp\sum_{\ell=0}^{\infty} \cbin{n+\ell-1}{\ell} t^\ell\rp\\
&= \sum_{N=1}^{\infty} \frac{t^N}{N} \lp H_0 + 
\sum_{r=1}^N (-1)^r N \cbin{N-1}{r-1} \frac{H_r}{[r]} \rp \\
&= \sum_{N=1}^{\infty} \frac{t^N}{N} \lp H_0 + 
 \sum_{r=1}^N (-1)^r  \cbin{N}{r} \frac{r}{[r]} H_r \rp
= \sum_{N=1}^{\infty} \frac{t^N}{N} \wt{H}_N
\end{align*}

and \eqref{eq:euler-pf} follows.

\end{pf}

\section{Proof of Theorem \ref{thm:cp-series} I}
\label{sec:proof}

In this section we prove Theorem \ref{thm:cp-series}
\eqref{it:rationality} and \eqref{it:limits}
for $\g=\Lsl_2$.
Our main tool is a straightening identity given in Proposition \ref{pr:straight}.
Denote the loop generators of $\qloopsl{2}$ by $\{E_k,F_k,H_k\}_{k\in\Z}$,
and set $K = q^{H_0}$.

There is an isomorphism of algebras $\Omega : \qloopsl{2} \to U_{q^{-1}}(L\Lsl_2)^{\op}$,
given on the loop generators by $\Omega(X_k) = X_{-k}$ ($X=E,F,H$)
(see \cite[Prop. 1.3]{cp-rootofunity}). Using this,
one obtains the $-$ case of
Theorem \ref{thm:cp-series} \eqref{it:rationality} from its $+$ counterpart.

\subsection{A commutation relation}\label{ssec:comm}
We begin by reviewing the commutation relation between
$\PP^+(z)$ and the raising/lowering
operators  of $\qloopsl{2}$ 
obtained in \cite[Lemma 3.3]{cp-rootofunity}.

\begin{lem}\label{lem:comm}
Let $E(w) = \sum_{n\in\Z} E_nz^{-n}$ and
$F(w) = \sum_{n\in\Z} F_nz^{-n}$. Then, 
\begin{align*}
\Ad(\PP^+(z))\cdot E(w) &= (1-q^2wz^{-1})(1-wz^{-1}) E(w)\\
\Ad(\PP^+(z))^{-1}\cdot F(w) &= (1-q^2wz^{-1})(1-wz^{-1}) F(w)
\end{align*}
\end{lem}

\begin{pf}
Using the commutation relation $\Ad(\cg^+(z))\cdot E(w) = \frac{q^2z-w}{z-q^2w} E(w)$,
and letting $\cgn^+(z) = q^{-H_0}\cg^+(z)$, we get:
\begin{align*}
\Ad(\cgn^+(z))\cdot E(w) &= \frac{1-q^{-2}wz^{-1}}{1-q^2wz^{-1}} E(w) 
= \frac{(1-q^{-2}wz^{-1})(1-wz^{-1})}{(1-q^2wz^{-1})(1-wz^{-1})} E(w) \\
&= \frac{p(q^2z,w)}{p(z,w)} E(w)\,
\end{align*}
where $p(z,w) = (1-q^2wz^{-1})(1-wz^{-1})$. Comparing with
the difference equation defining $\PP^+(z)$, we conclude that
$\Ad(\PP^+(z))\cdot E(w) = p(z,w) E(w)$. The argument for $F's$
is similar.
\end{pf}

In the proof of Proposition \ref{pr:straight} below, only the $F$ case
of the following identites is needed.

\begin{cor}\label{cor:comm}
\[
\begin{aligned}
\Ad(\PP^+(z)^{-1})\cdot \lp\sum_{n=0}^{\infty} F_{n+1}z^{-n}\rp
&= F_1 - q^2 F_2z^{-1}\\
\Ad(\PP^+(z))\cdot \lp\sum_{n=0}^{\infty} E_{n+1}z^{-n}\rp
&= E_1 - q^2 E_2z^{-1}
\end{aligned}
\]
\end{cor}

\begin{pf}
By the lemma above
\[
\Ad(\PP^+(z)^{-1})\cdot F_n = F_n - (1+q^2)F_{n+1}z^{-1} + q^2F_{n+2}z^{-2}\,,
\text{ for any } n\in\Z
\]
Therefore, we get
\begin{gather*}
\Ad(\PP^+(z)^{-1})\cdot (\sum_{n=0}^{\infty} F_{n+1}z^{-n}) = \\ 
F_1 + z^{-1}F_2 (1-(1+q^2)) + \sum_{m\geq 3} F_mz^{-m+1}
(1-(1+q^2)+q^2)
\end{gather*}
The proof of the second equation is verbatim, hence omitted.
\end{pf}

\subsection{A straightening formula}\label{ssec:straight}
Now we can state and prove the fundamental identity in
$\power{\qloopsl{2}}{z^{-1}}$ using
Chari--Pressley's straightening formula
\cite[Lemma 5.1]{cp-rootofunity}.
As before, $E^+(z) = \sum_{n=0}^{\infty} E_n z^{-n}$ and
we set $\xi^{(r)} = \frac{\xi^r}{[r]!}$.

\begin{prop}\label{pr:straight}
The following identity holds in $\qloopsl{2}[\![z^{-1}]\!]$.
\[
\begin{aligned}
\sum_{n=0}^{\infty} (-1)^n q^{n^2} E_0^{(n)}F_1^{(n)} K^{-n}z^{-n} = \hspace*{2.5in} \\
\hspace*{0.5in} \sum_{\ell=0}^{\infty} (-1)^{\ell} q^{\ell^2}
K^{-\ell} z^{-\ell} \cdot 
\PP^+(q^{-2\ell}z) \lp F_1-q^{2\ell+2}F_2 z^{-1}\rp^{(\ell)}
E^+(q^{-2\ell}z)^{(\ell)}   
\end{aligned}
\]
\end{prop}

\begin{pf}
We recall the straightening formula \cite[Lemma 5.1]{cp-rootofunity}.
The notations are explained below.
\begin{equation}\label{eq:cp-basic}
E_0^{(r)}F_1^{(s)} = \sum_{t=0}^{\text{min}(r,s)} (-1)^t q^{-t(r+s-t)}
\sum_{\begin{subarray}{c} a+b=t \\ a,b\in\N\end{subarray}}
D_a^-(\xi^{(s-t)}) K^t \mathbb{D}_b^+(\xi^{(r-t)})
\end{equation}

Following the conventions of \cite[\S 4]{cp-rootofunity},
we let $\xi$ be an indeterminate, and consider the algebra
homomorphisms $D^{\pm}(z) : \C[\xi] \to \power{\qloopsl{2}}{z^{-1}}$
defined by:
\[
D^+(z) \cdot \xi = \sum_{n=0}^{\infty} E_n z^{-n} = E^+(z) \quad\text{and}\quad
D^-(z) \cdot \xi = \sum_{n=0}^{\infty} F_{n+1}z^{-n} =: F_\infty(z)
\]
Note that $F_{\infty}(z) = z(F^+(z)-F_0)$ is the series appearing in
Corollary \ref{cor:comm} above.
Let $\mathbb{D}^+(z) := \LM(\PP^+(z)) D^+(z)$, where $\LM(X)$
is the operator of left multiplication by $X$. The operators $D_n^-$
and $\mathbb{D}_m^+$ appearing in \eqref{eq:cp-basic} above,
are the coefficients: $D^-(z) = \sum_{n=0}^{\infty} D_n^- z^{-n}$
and $\mathbb{D}^+(z) = \sum_{n=0}^{\infty} \mathbb{D}_n^+ z^{-n}$.\\

Note that $\Ad(K)\cdot (D^{\pm}(z)(\xi)) = q^{\pm 2}D^{\pm}(z)(\xi)$.
Since $D^{\pm}(z)$ are algebra homomorphisms, this implies that
$\Ad(K)\cdot (D^{\pm}(z)(\xi^r)) = q^{\pm 2r} D^{\pm}(z)(\xi^r)$, for
each $r\in\Z_{\geq 0}$. Equating coefficients of $z^{-n}$ on both sides
of this relation, one obtains the same for $D_n^{\pm}(\xi^r)$. Using this
observation, the fact that $K$ commutes with $\PP^+(z)$, and
\eqref{eq:cp-basic}, we are able to carry out the following computation:
\begin{align*}
(-1)^n q^{n^2} E_0^{(n)} F_1^{(n)} K^{-n} &=
\sum_{t=0}^n (-1)^{n-t} q^{(n-t)^2} \sum_{a+b=t}
D_a^-(\xi^{(n-t)}) K^t \mathbb{D}_b^+(\xi^{(n-t)}) K^{-n} \\
&= \sum_{t=0}^n (-1)^{n-t} q^{n^2-t^2} K^{-n+t} \sum_{a+b=t}
D_a^-(\xi^{(n-t)}) \mathbb{D}_b^+(\xi^{(n-t)})
\end{align*}
Changing $\ell = n-t$, we can
rewrite this as:
\[
(-1)^n q^{n^2} E_0^{(n)} F_1^{(n)} K^{-n} =
\sum_{\ell=0}^n (-1)^{\ell} q^{\ell(2n-\ell)} K^{-\ell}
\sum_{a+b=n-\ell} D_a^-(\xi^{(\ell)})\mathbb{D}^+_b(\xi^{(\ell)})
\]
Now we consider the series
\begin{align*}
\mathcal{E}(z) &:= \sum_{n=0}^{\infty} (-1)^n q^{n^2} E_0^{(n)}
F_1^{(n)} K^{-n} z^{-n} \\
&= \sum_{n=0}^{\infty} z^{-n} \cdot
\lp
\sum_{\ell=0}^n (-1)^{\ell} q^{\ell(2n-\ell)} K^{-\ell}
\sum_{a+b=n-\ell} D_a^-(\xi^{(\ell)})\mathbb{D}^+_b(\xi^{(\ell)})
\rp \\
&= \sum_{n=0}^{\infty} \lp \sum_{\ell=0}^n z^{-\ell} (-1)^{\ell}
q^{\ell^2} K^{-\ell} \sum_{a+b=n-\ell}
q^{2\ell a} z^{-a} D_a^-(\xi^{(\ell)}) q^{2\ell b}z^{-b}
\mathbb{D}^+_b(\xi^{(\ell)}) \rp \\
&= \sum_{\ell=0}^{\infty} (-1)^{\ell} q^{\ell^2} K^{-\ell} z^{-\ell}
\cdot \lp
D^-(q^{-2\ell}z)(\xi)^{(\ell)} \PP^+(q^{-2\ell}z) D^+(q^{-2\ell}z)(\xi)^{(\ell)}
\rp \\
&= \sum_{\ell=0}^{\infty} (-1)^{\ell} q^{\ell^2} K^{-\ell} z^{-\ell}
\cdot \lp
F_\infty\lp q^{-2\ell}z\rp^{(\ell)} \PP^+(q^{-2\ell}z) E^+\lp q^{-2\ell}z\rp^{(\ell)}
\rp 
\end{align*}
where we have used the fact that $D^{\pm}(w)$ are algebra homomorphisms,
and $\mathbb{D}^+(w)(f) = \PP^+(w)D^+(w)(f)$, for every $f\in\C[\xi]$.
The claim in the lemma now follows from Corollary \ref{cor:comm}.
\end{pf}

\subsection{Proof of Theorem \ref{thm:cp-series}, \eqref{it:rationality}}
\label{ssec:pf-1}

Let $V$ be a \fd type I representation of $\qloopsl{2}$, and
consider the evaluation of the identity from Proposition \ref{pr:straight}
on $V$. By the usual weight reasons, both sides of the identity become
finite sums. The left--hand side is clearly a polynomial in $z^{-1}$,
with coefficients from $\End(V)$. We write the right--hand side, truncated
at $M$ (depending on $V$):
\[
\PP^+(z)\cdot \lp \sum_{\ell=0}^M (-1)^{\ell} q^{\ell^2}K^{-\ell}z^{-\ell}
\PP^+(z)^{-1}\PP^+(q^{-2\ell} z) (F_1-q^2F_2z^{-1})^{(\ell)}E^+(q^{-2\ell}z)^{(\ell)}
\rp
\]
We claim that each term in the summation above is a rational function of
$z$, once evaluated on $V$. For $E^+(z)$ 
this holds by Proposition \ref{pr:rat currents}. 
Using
the difference equation $\PP^+(q^2z)\PP^+(z)^{-1} = \ol{\cg^+}(z) = q^{-H_0}\cg^+(z)$,
we have
\[
\PP^+(z)^{-1}\PP^+(q^{-2\ell}z) = \ol{\cg^+}(q^{-2}z)\cdots \ol{\cg^+}(q^{-2\ell}z)
\]
which is a rational function of $z$ by Proposition \ref{pr:rat currents}.
Therefore, we obtain an equation of the form $A(z) = \PP^+(z)B(z)$, where
both $A(z)$ and $B(z)$ are rational $\End(V)$--valued functions, taking
value $\Id_V$ at $z=\infty$, in particular, generically invertible.
Hence, $\PP^+(z)=A(z)B(z)^{-1}$ evaluated on
$V$, is a ``ratio" of two rational $\End(V)$--valued functions, and thus
is itself rational.

\subsection{Proof of Theorem \ref{thm:cp-series}, \eqref{it:limits}}
\label{ssec:pf-2}
Both $\PP^+(z)$ and $\PP^-(z)$ are viewed as rational
$\End(V)$--valued functions, taking value $\Id_V$ at $z=\infty,0$
respectively. Moreover, they satisfy related difference equations:
\[
\PP^{\pm}(q^2z) = q^{\mp H_0}\cg^\pm(z)\PP^\pm(z)
\]

By Proposition \ref{pr:rat currents} $\cg^+(z) = \cg^-(z)$
as rational $\End(V)$--valued functions, so that
\[
\PP^+(q^2z)\PP^-(q^2z)^{-1} = q^{-2H_0} \PP^+(z)\PP^-(z)^{-1}
\]

Thus $\cplim(z) := z^{H_0}\PP^+(z)\PP^-(z)^{-1}$ is a $q^2$--periodic,
rational $\End(V)$--valued function of $z$, hence it must be a constant:
$\cplim(z)=\cplim\in\End(V)$, for every $z\in\C$. 
To see that it is invertible, we
observe that, being rational functions invertible at $\infty,0$ respectively, 
$\PP^{\pm}(z)$ are invertible at all but finitely many values of $z\in\C$. 
The value of $\cplim(z)$ can be computed at a generic $z_0\in\C$ where
all the operators appearing in its definition are invertible.

Finally, we note that $z^{H_0}\PP^+(z) = \cplim\PP^-(z)$. The right--hand side
takes value $\cplim$ at $z=0$, so we get the existence of the following limit:
\[
\lim_{z\to 0} z^{H_0}\PP^+(z) = \cplim
\]

Similarly, $\cplim^{-1} = \lim_{z\to\infty} z^{-H_0} \PP^-(z)$.

\section{Proof of Theorem \ref{thm:cp-series} II}
\label{ssec:pf-3}

Let $\{\E_i,\F_i,\K_i\}_{i=0,1}$ denote the Chevalley generators
of $\qloopsl{2}$.
The following isomorphism will be used in our proof of
Theorem \ref{thm:cp-series} \eqref{it:lattice}. 
We recall that this isomorphism is
$\biso{-}$ in the notations of Section \ref{ssec:beck}.
\[
\begin{aligned}
\E_1 &= E_0 & \F_1 &= F_0 & \H_1 &= H_0 \\
\E_0 &= K^{-1}F_1 & \F_0 &= E_{-1}K & \H_0 &= -H_0
\end{aligned}
\]

Recall that we have to prove the following identity, where both
sides are viewed as operators on $V\in\Rloop$, see Example \ref{ex:sl2}
and equation \eqref{eq: triple-ex} in Section \ref{ssec: q-weyl-group}.

\begin{equation}\label{eq:pf-3-1}
\qS_1^{-1}\qS_0^{-1} = (-q)^{-H_0}\lim_{z\to 0} z^{H_0}\PP^+(z)
=(-q)^{-H_0}\cplim
\end{equation}
where, as before $\cplim=\lim_{z\to 0}z^{H_0}\PP^+(z)$.

\subsection{Outline of the proof}
Our proof of \eqref{eq:pf-3-1} is in the following steps. 
First, we show in Section \ref{ssec:pf-3-comm} that both sides
satisfy the same commutation relation with the loop generators of
$\qloopsl{2}$. Thus, it suffices
to prove \eqref{eq:pf-3-1} on a subspace of $V$ which generates
it as $\qloopsl{2}$--representation.

Next, we take $V' := \Ker(E_{-1}) = \Ker(\F_0) \subset V$ as the
generating subspace, and
focus on a typical weight component $V'[\ws]=V'\cap V[\ws]$,
where
$V[\ws]=\{v\in V : H_0\cdot v = \ws v\}$.
Note that, by finite--dimensionality of $V$, $\ws\in\N$.
Using the fact that $E_0^{(\ws)}$ is invertible on $V[-\ws]$,
equation \eqref{eq:pf-3-1} on $V'[\ws]$ is implied by the
following two, proved in Sections \ref{ssec:pf-3I} and \ref{ssec:pf-3II}
respectively.

\[
E_0^{(\ws)}\qS_0^{-1} = \lp\sum_{\ell=0}^{\infty}
(-1)^{\ell} q^{\ell(\ws+1)} F_0^{(\ell)}E_0^{(\ell)}\rp \cplim,\ \ 
E_0^{(\ws)}\qS_1 = (-q)^{\ws} \sum_{\ell=0}^{\infty}
(-1)^{\ell} q^{\ell(\ws+1)} F_0^{(\ell)}E_0^{(\ell)}
\]

\subsection{Commutation relations}\label{ssec:pf-3-comm}

It follows from Lemma \ref{lem:comm} that
\begin{equation}\label{eq:C-comm}
\Ad(\cplim)\cdot E_k = q^2E_{k+2}\ \ \  \Ad(\cplim)\cdot F_k = q^{-2}F_{k-2}
\end{equation}
It is clear that $\cplim$ commutes with $H_k$, for every $k\in\Z$.\\

The computation of $\Ad(\qS_0\qS_1)$ acting on $\qloopsl{2}$ can
be found in \cite[Lemmas 9.2, 9.3, \S 9.4]{sachin-valerio-3}:

\begin{equation}\label{eq:L-comm}
	\Ad(\qS_0\qS_1)\cdot E_k = E_{k-2}\,\quad
	\Ad(\qS_0\qS_1)\cdot F_k = F_{k+2}\, \quad
	\Ad(\qS_0\qS_1)\cdot H_k = H_k
\end{equation}

Hence, our claim that $\Ad(\qS_1^{-1}\qS_0^{-1})=\Ad((-q)^{-H_0}\cplim)$
on $\qloopsl{2}$ follows.

\subsection{$U_q(\Lsl_2)$ computation}\label{ssec:pf-3-sl2}
The following computation is standard and will be needed in our proof
(see, for instance, \cite[\S 8.3]{jantzen}).

Let $\IR_n$ be the $(n+1)$--dimensional, irreducible
representation of $U_q(\sl_2)$. There is a basis
$\{m_n(r) : 0\leq r\leq n\}$ of $\IR_n$, in which the $U_q(\sl_2)$--action
is given by:

\[
\K m_n(r) = q^{n-2r} m_n(r)\ \ \ 
\E m_n(r) = [n-r+1] m_n(r-1)\ \ \ 
\F m_n(r) = [r+1] m_n(r+1)
\]

\begin{prop}\label{pr: sl2}
With the notations as above, we have
\[
\qS \cdot m_n(r) = (-1)^{n-r} q^{(n-r)(r+1)} m_n(n-r)
\]
\end{prop}

\subsection{Limit of the straightening identity}
\label{ssec:pf-3I}
Fix $\ws\in\N$ and recall that
\[
V'[\ws]=\{v\in V : H_0\cdot v = \ws v\text{ and }
\F_0\cdot v = 0\}
\]

\begin{lem}\label{lem:C-straight}
We have the following identity on $V'[\ws]$:
\begin{equation}\label{eq:pf-3-2}
E_0^{(\ws)}\qS_0^{-1} = \lp\sum_{\ell=0}^{\infty}
(-1)^{\ell} q^{\ell(\ws+1)} F_0^{(\ell)}E_0^{(\ell)}\rp \cplim
\end{equation}
\end{lem}

\begin{pf}
Consider the equation obtained in Proposition \ref{pr:straight}:
\begin{equation}\label{eq:straight}
\begin{aligned}
\sum_{n=0}^{\infty} (-1)^n q^{n^2} E_0^{(n)}F_1^{(n)} K^{-n}z^{-n} = \hspace*{2.5in} \\
\hspace*{0.5in} \sum_{\ell=0}^{\infty} (-1)^{\ell} q^{\ell^2}
K^{-\ell} z^{-\ell} \cdot 
\PP^+(q^{-2\ell}z) \lp F_1-q^{2\ell+2}F_2 z^{-1}\rp^{(\ell)}
E^+(q^{-2\ell}z)^{(\ell)}   
\end{aligned}
\end{equation}

We evaluate both sides on $V'[\ws]$. Note that,
on this subspace, we have $\E_0^{(\ws+r)}=0$ for every
$r\geq 1$. 
Using, $\E_0 = K^{-1}F_1$, we get $\E_0^{(n)} = q^{n(n+1)}
F_1^{(n)}K^{-n}$, which shows that the left--hand side truncates
at $n=\ws$.
Thus, the limit $z\to 0$ of the left--hand side
of \eqref{eq:straight}, multiplied by
$z^{\ws}$ is:
\[
\lim_{z\to 0} (z^\ws \mathrm{L.H.S.}(z)) = 
(-1)^{\ws} q^{\ws^2} E_0^{(\ws)}F_1^{(\ws)} K^{-\ws}
\]
Combining this with the $U_q(\Lsl_2)$ fact,
directly following from Proposition \ref{pr: sl2}, that
$\qS_0^{-1} = (-1)^{\ws} q^{-\ws} \E_0^{(\ws)}$ on $\Ker(\F_0)
\cap V[\ws]$, we get

\begin{equation}\label{eq:pf-3I-1}
\lim_{z\to 0} (z^\ws \mathrm{L.H.S.}(z)) = 
E_0^{(\ws)}\qS_0^{-1}
\end{equation}

Consider now the right--hand side of \eqref{eq:straight}. 
First, we need the following equality
of rational functions (Proposition \ref{pr:rat currents})
\[
E^+(z) = \sum_{n=0}^{\infty} E_nz^{-n} = - \sum_{m=1}^{\infty} E_{-m} z^m
= E^-(z)
\]
Let $\ul{E^-}(z) = E^-(z)+E_{-1}z = -\sum_{m\geq 2} E_{-m}z^m$. Next, we need
the following commutation relation.\\

\noindent {\bf Claim.} For every $\ell\geq 0$, we have
\[
E^-(z)^{\ell} = q^{\ell(\ell-1)} \ul{E^-}(z)^{\ell} \text{ modulo
the left ideal generated by } E_{-1}\ .
\]

Given this claim, and replacing $K$ by $q^{\ws}$,
a typical summand of the right--hand side of \eqref{eq:straight} can be rewritten
as
\[
(-1)^{\ell} q^{\ell^2-\ell\ws+\ell(\ell-1)} \PP^+(q^{-2\ell}z) 
(zF_1-q^{2\ell+2}F_2)^{(\ell)}
\ul{E^-}(q^{-2\ell}z)^{(\ell)}\ z^{-2\ell}
\]
Note that $\ul{E^-}(q^{-2\ell}z)^{(\ell)} = (q^{-2\ell}z)^{2\ell} (-E_{-2}-
\sum_{r\geq 1} E_{-2-r}(q^{-2\ell}z)^r)^{(\ell)}$. 
Upon simplifying the exponent of $q$,
the $\ell$--th summand becomes
\[
(-1)^{\ell} q^{\ell-\ell\ws} \PP^+(q^{-2\ell}z)(F_2-q^{-2\ell-2}zF_1)^{(\ell)}
(E_{-2}+\sum_{r\geq 1} E_{-2-r}(q^{-2\ell}z)^r)^{(\ell)}
\]

With the exception of $\PP^+(q^{-2\ell}z)$ term, we can set $z=0$ in the
rational functions appearing above. Now we use
$$\lim_{z\to 0} z^{\ws} \PP^+(q^{-2\ell}z) = q^{2\ell\ws} 
\lim_{z\to 0} (q^{-2\ell}z)^{\ws} \PP^+(q^{-2\ell}z)
= q^{2\ell\ws}\cplim$$ 
to get the following answer for the right--hand side
of \eqref{eq:straight}.

\begin{align*}
\lim_{z\to 0} z^{\ws}\mathrm{R.H.S.}(z) &=
\cplim \cdot \lp \sum_{\ell=0}^{\infty} (-1)^{\ell} q^{\ell(\ws+1)}
F_2^{(\ell)} E_{-2}^{(\ell)}\rp \\
&= \lp \sum_{\ell=0}^{\infty} (-1)^{\ell} q^{\ell(\ws+1)}
F_0^{(\ell)}E_0^{(\ell)} \rp \cdot \cplim
\end{align*}
where, we used the commutation relation between $\cplim$ and $E_k,F_k$,
equation \eqref{eq:C-comm} above. Note that we have interchanged the order of
$\lim_{z\to 0}$ and $\sum_{\ell=0}^{\infty}$, which is justified since this
sum, evaluated on $V[\ws]$ is in fact finite.
The lemma is proved, modulo the claim.\\

\noindent {\em Proof of the claim.} Consider the commutation relation
\[
E_{-k}E_{-\ell-1} - q^2E_{-\ell-1}E_{-k} = 
q^2E_{-k-1}E_{-\ell} - E_{-\ell}E_{-k-1}
\]
Multiply both of its sides by $z^{k+1}w^{\ell+1}$ and sum over all
$k,\ell\geq 1$ to get:
\[
zE^-(z)\ul{E^-}(w) - q^2z\ul{E^-}(w)E^-(z) = 
q^2w\ul{E^-}(z)E^-(w) - wE^-(w)\ul{E^-}(z)\ .
\]
Now set $z=w$ to get $E^-(z)\ul{E^-}(z) = q^2 \ul{E^-}(z)E^-(z)$.
The claim follows by an easy induction on $\ell$
together with the fact that $E^-(z) = \ul{E^-}(z)$ modulo
the left ideal generated by $E_{-1}$.

\end{pf}

\subsection{End of the proof}
\label{ssec:pf-3II}
With the help of Lemma \ref{lem:C-straight}, the 
equation \eqref{eq:pf-3-1} follows from
the following computation of the action of $\qS_1$ on a positive weight
space $V[\ws]$ (this is purely a $U_q(\Lsl_2)$ statement),
and the fact that $E_0^{(\ws)}$ is an invertible operator on
the weight space $V[-\ws]$:

\begin{equation}\label{eq:pf-3-3}
E_0^{(\ws)}\qS_1 = (-q)^{\ws} \sum_{\ell=0}^{\infty}
(-1)^{\ell} q^{\ell(\ws+1)} F_0^{(\ell)}E_0^{(\ell)}
\end{equation}

\begin{pf}

Since $V$ is semisimple as $U_q(\Lsl_2)$ representation,
it is enough to check \eqref{eq:pf-3-3} on an irreducible
representation $\IR_{\Lambda}$, where $\Lambda\in\N$ (see
the notations before Proposition \ref{pr: sl2}), a positive
weight $\ws = \Lambda-2r$, where $0\leq r\leq \lfloor\frac{\Lambda}{2}\rfloor$.
Using the $U_q(\Lsl_2)$--action on a basis $\{m(r):0\leq r\leq \Lambda\}$,
we get (we have used $\Lambda-r = \ws+r$ here)
\[
E_0^{(\ws)}\qS_1 \cdot m(r) = 
(-1)^{\ws+r} q^{(\ws+r)(r+1)}\qbin{\ws+r}{r}{q} m(r)
\]

Similarly, the right--hand side of \eqref{eq:pf-3-3} evaluated
on $m(r)$ becomes
\[
(-1)^{\ws} q^{\ws} \sum_{\ell=0}^r
(-1)^{\ell} q^{\ell(\ws+1)} \qbin{r}{\ell}{q}
\qbin{\ws+r+\ell}{\ell}{q}
\]
Hence, equation \eqref{eq:pf-3-3} becomes the following
$q$--binomial identity (we set $y=\ws+r$)
\[
\sum_{\ell=0}^r (-1)^{\ell} q^{\ell(y-r+1)}
\qbin{r}{\ell}{q} \qbin{y+\ell}{\ell}{q} = 
(-1)^r q^{r(y+1)}\qbin{y}{r}{q}
\]

One way to obtain this equation, is to use the following
``iterated $q$--Pascal" relation (see, for instance,
\cite[\S 1.3.1, (e)]{lusztig-book}).

\[
\qbin{y+\ell}{\ell}{q} = \sum_{j=0}^{\ell}
q^{-(\ell-j)y+\ell j} \qbin{\ell}{j}{q}
\qbin{y}{j}{q}
\]

Substituting it in the left--hand side, and comparing
coefficients of $\ds\qbin{y}{j}{q}$, the equation reduces
to the following:
\begin{gather*}
(-1)^j \qbin{r}{j}{q} q^{j(y-r+j+1)}
\sum_{i=0}^{r-j} (-1)^i q^{-i(r-j-1)}\qbin{r-j}{i}{q}
 = \\
\left\{
\begin{array}{cl}
0 & \text{ if } 0\leq j<r \\
(-1)^r q^{r(y+1)} & \text{ if } j=r
\end{array}
\right.
\end{gather*}

This is again well--known (see, for instance, 
\cite[\S 1.3.4]{lusztig-book}), and this completes
our proof. 
\end{pf}

\section{Nakajima quiver varieties}\label{sec:NQV}

In this section we compute the action of the lattice operators on the equivariant
$K$--theory of Nakajima quiver varieties, assuming $\g$ is of type $\mathsf{ADE}$.

\subsection{Quiver varieties}\label{ssec:qv-defn}

We briefly review the salient features of quiver varieties, and refer to \cite{ginzburg-qv}
for more details.

Let $\Gamma=(\bfI,\bfE)$ be a quiver whose underlying unoriented graph is the
Dynkin diagram of $\g$, and $s,t:\bfE\to\bfI$ the source and target maps.
Let $\mathsf{v},\mathsf{w}\in\N^{\bfI}$ be dimension vectors, and $(V_k)_{k\in
\bfI}$ and $(W_k)_{k\in\bfI}$ two collections of vector spaces over $\C$ such
that $\dim(V_k)=\mathsf{v}_k$ and $\dim(W_k)=\mathsf{w}_k$. 
The quiver variety $\M(\sfv,\sfw)$ is the $\mathrm{GIT}$ quotient
of the Hamiltonian
$\GL(\sfv)=\prod_{k\in \bfI}\GL(\sfv_k)$--action on
the following symplectic vector space
\[\begin{split}
\mathbb{M}(\sfv,\sfw) = 
&\bigoplus_{a\in\bfE}
\Hom( V_{s(a)},V_{t(a)}) \oplus \Hom(V_{s(a^*)},V_{t(a^*)})\\
&\bigoplus_{k\in\bfI} \Hom(V_k,W_k)\oplus \Hom(W_k,V_k)
\end{split}\]
Here, for each $a\in\bfE$, $a^*$ is included as a new edge with the opposite
orientation of $a$, so $s(a^*)=t(a)$ and $t(a^*)=s(a)$.

A point of $\mathbb{M}(\sfv,\sfw)$ is denoted by $(B,i,j)$ where,
\begin{itemize}
\item $B = (B_a:V_{s(a)}\to V_{t(a)})_{a\in\bfE\cup\bfE^*}$
\item $i = (i_k:W_k\to V_k)_{k\in\bfI}$
\item $j = (j_k:V_k\to W_k)_{k\in\bfI}$
\end{itemize}
Thus, a point in $\M(\sfv,\sfw)$ is the $\GL(\sfv)$--orbit
of $[(B,i,j)]$, where $(B,i,j)$ satisfies:
\begin{itemize}
\item ADHM equations. $\mu(B,i,j)=0$, where
the moment map $\mu:\mathbb{M}(\sfv,\sfw)\to\mathfrak{gl}(\sfv)$
is given by:
\[
\mu(B,i,j) = \lp
\sum_{\begin{subarray}{c} a\in\bfE \\ t(a)=k\end{subarray}}
B_aB_{a^*} - 
\sum_{\begin{subarray}{c} a\in\bfE \\ s(a)=k\end{subarray}}
B_{a^*}B_a + 
i_kj_k
\rp_{k\in\bfI}
\]
\item $(B,i,j)\in\mathbb{M}(\sfv,\sfw)$ is {\em stable}
\cite[Defn.~3.2.1]{nakajima-qaffine}, that is,
every $B$--invariant subspace of $(V_k)$, which is contained
in the kernel of $j$, is zero.
\end{itemize}

Let $\mu^{-1}(0)^s\subset\mathbb{M}(\sfv,\sfw)$ denote the set of
all points satisfying the two conditions written above. Then,
$\M(\sfv,\sfw)=\mu^{-1}(0)^s/\GL(\sfv)$.
We summarize its key properties as follows, see \cite[Prop.~2.3.2]{nakajima-qaffine}.

\begin{itemize}
\item $\M(\sfv,\sfw)$ is a connected, quasi--projective, smooth,
symplectic variety. 

\item $\M(\sfv,\sfw)$ comes equipped with an action of
$\GL(\sfw)\times \nC$ (see \cite[\S2.7]{nakajima-qaffine}),
where $z\in\nC$ acts by scaling all the linear maps
\[z\cdot [(B,i,j)] = [(z B, z i, z j)]\]

\item $\mu^{-1}(0)^s \to \M(\sfv,\sfw)$ is a principal $\GL(\sfv)$--bundle.
\end{itemize}

\subsection{Tautological bundles}\label{ssec:bundles}
For each $k\in \bfI$, there is a natural action of
$\GL(\sfv)$ on $V_k$. Since $\mu^{-1}(0)^s \to \M$ is
a principal $\GL(\sfv)$--bundle, we can define a vector
bundle $\V_k$ on $\M(\sfv,\sfw)$ by
\[
\V_k := \mu^{-1}(0)^s \times_{\GL(\sfv)} V_k
\]
which is $\GL(\sfw)\times \nC$--equivariant 
where $\GL(\sfw)$ acts trivially on the fiber.
Similarly, we define $\W_k$ to be the trivial
vector bundle $\M(\sfv,\sfw)\times W_k$ together with the
$\GL(\sfw)\times \nC$--equivariant structure coming from
the natural action
of $\GL(\sfw)$ on $W_k$ and trivial $\nC$ action on the
fiber.

For any $m\in \Z$, let $L(m)$ be the one--dimensional representation of 
$\nC$ given by $z\to z^m$ and, for any $\nC$--module $A$, set $\ds q^
m A = L(m)\otimes A$. Let $C_k(\sfv,\sfw)$ be the complex of equivariant
vector bundles on $\M(\sfv,\sfw)$ given by \cite[(2.9.1)]{nakajima-qaffine}
\[C_k(\sfv,\sfw):\qquad
q^{-2}\V_k \stackrel{\sigma_k}{\longrightarrow}
q^{-1}\lp\W_k \oplus \bigoplus_{\ell:a_{k\ell}=-1} \V_\ell \rp
\stackrel{\tau_k}{\longrightarrow}
\V_k\]
where $(a_{k\ell})_{k,\ell\in\bfI}$ is the Cartan matrix of $\g$, and the terms
are in degrees $-1,0,1$. The maps $\sigma_k,\tau_k$ on the fiber over a
point $[(B,i,j)]\in\M(\sfv,\sfw)$ are given by
\[\sigma_k =  j_k+
\bigoplus_{\begin{subarray}{c} a\in\bfE\cup\bfE^* \\
t(a)=k\end{subarray}} B_{a^*}
\aand
\tau_k =  i_k+\sum_{\begin{subarray}{c} a\in\bfE\cup\bfE^* \\
s(a)=k\end{subarray}} \varepsilon(a)B_a
\]
where $\varepsilon(a)=1$ if $a\in\bfE$, and $-1$ if $a\in\bfE^*$.
Note that 
\[
\rk(C_k(\sfv,\sfw)) = \sfw_k-2\sfv_k+
\sum_{\ell:a_{k\ell}=-1} \sfv_\ell = \langle \Lambda_{\sfw}-\alpha_{\sfv},\alpha_k\rangle
\]
where $\Lambda_{\sfw} = \sum_{i\in\bfI} \sfw_i\varpi_i$ and
$\alpha_{\sfv}=\sum_{i\in\bfI} \sfv_i\alpha_i$.

\subsection{Action of $\qloopo$}\label{ssec:U0-action}

Recall that if $X$ is a smooth, quasi--projective variety over $\C$ endowed with
an algebraic action of a linear algebraic group $G$, 
$K_G(X)$ denotes the Grothendieck group of $G$--equivariant coherent sheaves
on $X$. $K_G(X)$ is a module over the representation ring $R(G) = K_G(\mathrm
{pt})$, and is generated 
by the equivalence classes of $G$--equivariant, algebraic vector bundles on $X$
(the reader should consult \cite[\S~6.1]{nakajima-qaffine} for a quick summary of
equivariant $K$--theory, and references for the subject).

In \cite[\S 9]{nakajima-qaffine}, Nakajima defines an action of $\qloop$ on the
equivariant $K$--theory of quiver varieties. More precisely, for a fixed $\sfw$,
there is an algebra homomorphism
\[\Nqv : \qloop \to \End\lp \bigoplus_{\sfv\in\N^\bfI} 
K_{\GL(\sfw)\times \nC}(\M(\sfv,\sfw))_\C\rp\]
where $K_G(X)_\C := K_G(X)\otimes_\Z \C$.
For our purposes, only the action of $\qloopo$
is relevant, and is
given as follows. Recall that $\qloopo$ is the commutative subalgebra
of $\qloop$ generated by $\{\psi_{k,n}^\pm\}_{k\in \bfI, n\in \Z_{\geq 0}}$.

Denote by
\[\textstyle{\bigwedge_u}:K_G(X)\to [\mathcal{O}_X]+K_G(X)\otimes u\Z[\![u]\!]\]
the group homomorphism given by $\bigwedge_u(E)=\sum_i u^i\bigwedge
^i E$ for any $G$--equivariant vector bundle $E$ on $X$. Then, for any
$k\in \bfI$, the generating series $\cg^{\pm}_k(z)$ acts on $K_{\GL(\sfw)
\times \nC}(\M(\sfv,\sfw))_\C$ as multiplication by the following element \cite
[(9.2.1)]{nakajima-qaffine}

\begin{equation}\label{eq: action-psi}
\cg^{\pm}_k(z) \mapsto q^{\langle\Lambda_{\sfw}-\alpha_{\sfv},\alpha_k\rangle}
\lp\frac{\bigwedge_{-1/qz} C_k(\sfv,\sfw)}
{\bigwedge_{-q/z} C_k(\sfv,\sfw)}\rp^{\pm}
\end{equation}

Note that the action of $\qloopo$ preserves the direct sum decomposition \ie
for any fixed $\sfv,\sfw$, $\Nqv$ maps $\qloopo$ to $\End\lp K_{\GL(\sfw)\times
\nC}(\M(\sfv,\sfw))_\C\rp$.

\subsection{Lattice operators}\label{ssec:nakajima-lattice}

The following result computes the \qWg action of the coroot lattice $Q^\vee$
on $K_{\GL(\sfw)\times\nC}(\M(\sfv,\sfw))_\C$. As in Theorem \ref{thm:cp-series},
this action depends on the chosen sign $\bsign$, and is denoted by
$\lambdao$ for notational simplicity.

\begin{thm}\label{thm:nakajima-lattice}
Let $\qloopo$ act on $K_{\GL(\sfw)\times\nC}(\M(\sfv,\sfw))_\C$ via the action
homomorphism $\Nqv$. Then, for every $k\in\bfI$, the following holds
\begin{enumerate}
\item The series $\PP_k^+(z)$ acts by tensoring with $\bigwedge_{-q/z} C_k(\sfv,\sfw)$.
\item The element $(-\bsign(k))^{H_{k,0}}\lambdao(L_k)$ acts by tensoring
with the determinant line bundle
\[\det(C_k(\sfv,\sfw))^* = 
q^{\langle\Lambda_{\sfw}-\alpha_{\sfv},\alpha_k\rangle} 
\det(\W_k)^*\otimes \det(\V_k)^{\otimes 2}\otimes
\lp\bigotimes_{\ell:a_{k\ell}=-1} \det(\V_\ell)^*\rp
\]
\end{enumerate}
\end{thm}
\begin{pf}
(1) Under \eqref{eq: action-psi}, the action of $\psi_{k,0}^+$ is given by $q^{\rk
C_k(\sfv,\sfw)}$. The claim then follows from the uniqueness of solution near
$\infty$ of the $q$--difference equation \eqref{eq:cp-diff}.

(2) If $E$ is a rank $r$ vector bundle, one has
\[\lim_{z\to 0}\,z^r\textstyle{\bigwedge_{-q/z}}(E) = (-1)^rq^r\bigwedge^r(E) = (-q)^r\det(E)\]
Thus, the automorphism $\ds\cplim_k=\lim_{z\to 0}z^{\HH_k}\PP_k^+(z)$ is
given by $(-q)^{\rk C_k(\sfv,\sfw)}\det(C_k(\sfv,\sfw))$. The result now follows
from Theorem \ref{thm:cp-series}.
\end{pf}


\begin{thebibliography}{EGH{\etalchar{+}}11}

\bibitem[Bec94]{beck-braid}
J.~Beck, \emph{Braid group action and quantum affine algebras}, Comm. Math.
  Phys. \textbf{165} (1994), no.~3, 555--568.

\bibitem[BK96]{beck-kac}
J.~Beck and V.~G. Kac, \emph{Finite-dimensional representations of quantum
  affine algebras at roots of unity}, J. Amer. Math. Soc. \textbf{9} (1996),
  no.~2, 391--423.

\bibitem[Bou02]{bourbaki-lie456}
N.~Bourbaki, \emph{Lie groups and {L}ie algebras. {C}hapters 4--6}, Elements of
  Mathematics (Berlin), Springer-Verlag, Berlin, 2002, Translated from the 1968
  French original by Andrew Pressley.

\bibitem[Cha86]{chari-integrable}
V.~Chari, \emph{Integrable representations of affine {L}ie-algebras}, Invent.
  Math. \textbf{85} (1986), no.~2, 317--335.

\bibitem[CKL13]{CKL}
S.~Cautis, J.~Kamnitzer, and A.~Licata, \emph{Coherent sheaves on quiver
  varieties and categorification}, Math. Ann. \textbf{357} (2013), no.~3,
  805--854.

\bibitem[CP91]{cp-qaffine}
V.~Chari and A.~Pressley, \emph{Quantum affine algebras}, Comm. Math. Phys.
  \textbf{142} (1991), no.~2, 261--283.

\bibitem[CP97]{cp-rootofunity}
\bysame, \emph{Quantum affine algebras at roots of unity}, Represent. Theory
  \textbf{1} (1997), 280--328.

\bibitem[CR08]{chuang-rouquier}
J.~Chuang and R.~Rouquier, \emph{Derived equivalences for symmetric groups and
  {$\mathfrak {sl}_2$}-categorification}, Ann. of Math. (2) \textbf{167}
  (2008), no.~1, 245--298.

\bibitem[Dri87]{drinfeld-yangian-qaffine}
V.~G. Drinfeld, \emph{A new realization of {Y}angians and of quantum affine
  algebras}, Dokl. Akad. Nauk SSSR \textbf{296} (1987), no.~1, 13--17.

\bibitem[EGH{\etalchar{+}}11]{pasha-book}
P.~Etingof, O.~Golberg, S.~Hensel, T.~Liu, A.~Schwendner, D.~Vaintrob, and
  E.~Yudovina, \emph{Introduction to representation theory}, Student
  Mathematical Library, vol.~59, American Mathematical Society, Providence, RI,
  2011, With historical interludes by Slava Gerovitch.

\bibitem[EV02]{etingof-varchenko}
P.~Etingof and A.~Varchenko, \emph{Dynamical {W}eyl groups and applications},
  Adv. Math. \textbf{167} (2002), no.~1, 74--127. 

\bibitem[FH15]{frenkel-hernandez}
E.~Frenkel and D.~Hernandez, \emph{Baxter's relations and spectra of quantum
  integrable models}, Duke Math. J. \textbf{164} (2015), 2407--2460.

\bibitem[FH25]{frenkel-hernandez-recent}
\bysame, \emph{Extremal monomial property of
  {$q$}-characters and polynomiality of the {$X$}-series}, J. Algebraic Combin.
  \textbf{62} (2025), Paper No. 41, 40. 

\bibitem[FR99]{frenkel-reshetikhin-qchar}
E.~Frenkel and N.~Yu. Reshetikhin, \emph{{The $q$-characters of representations
  of quantum affine algebras and deformed $\mathcal{W}$ algebras}}, {Recent
  Developments in Quantum Affine Algebras and Related Topics} (Raleigh, NC),
  Contemp. Math., vol. 248, AMS, Providence, RI, 1999, pp.~163--205.

\bibitem[Gar80]{garland}
H.~Garland, \emph{The arithmetic theory of loop groups}, Inst. Hautes
  \'{E}tudes Sci. Publ. Math. (1980), no.~52, 5--136.

\bibitem[Gin12]{ginzburg-qv}
V.~Ginzburg, \emph{Lectures on {N}akajima's quiver varieties}, Geometric
  methods in representation theory. {I}, S\'{e}min. Congr., vol. 24-I, Soc.
  Math. France, Paris, 2012, pp.~145--219.

\bibitem[GTL13]{sachin-valerio-3}
S.~Gautam and V.~Toledano~Laredo, \emph{Monodromy of the trigonometric
  {C}asimir connection for {$\mathfrak{sl}_2$}}, Noncommutative birational
  geometry, representations and combinatorics, Contemp. Math., vol. 592, Amer.
  Math. Soc., Providence, RI, 2013, pp.~137--176.

\bibitem[GTL16]{sachin-valerio-2}
\bysame, \emph{Yangians, quantum loop algebras, and abelian difference
  equations}, J. Amer. Math. Soc. \textbf{29} (2016), no.~3, 775--824.

\bibitem[GTL17]{GTL17}
\bysame, \emph{Meromorphic tensor equivalence for
  {Y}angians and quantum loop algebras}, Publ. Math. Inst. Hautes \'Etudes Sci.
  \textbf{125} (2017), 267--337. 

\bibitem[Har92]{hardy-divergent}
G.~H. Hardy, \emph{Divergent series}, \'{E}ditions Jacques Gabay, Sceaux, 1992,
  With a preface by J. E. Littlewood and a note by L. S. Bosanquet, Reprint of
  the revised (1963) edition.

\bibitem[Her05]{hernandez-05}
D.~Hernandez, \emph{Representations of quantum affinizations and fusion
  product}, Transform. Groups \textbf{10} (2005), no.~2, 163--200. 

\bibitem[Her07]{hernandez-drinfeld-coproduct}
\bysame, \emph{Drinfeld coproduct, quantum fusion tensor category and
  applications}, Proc. Lond. Math. Soc. (3) \textbf{95} (2007), no.~3,
  567--608.


\bibitem[Kac90]{kac}
V.~G. Kac, \emph{Infinite-dimensional {L}ie algebras}, third ed., Cambridge
  University Press, Cambridge, 1990.

\bibitem[KR90]{kirillov-reshetikhin}
A.~N. Kirillov and N.~Reshetikhin, \emph{{$q$}-{W}eyl group and a
  multiplicative formula for universal {$R$}-matrices}, Comm. Math. Phys.
  \textbf{134} (1990), no.~2, 421--431.

\bibitem[Lus88]{lusztig-qdef}
G.~Lusztig, \emph{Quantum deformations of certain simple modules over
  enveloping algebras}, Adv. in Math. \textbf{70} (1988), no.~2, 237--249.

\bibitem[Lus90]{lusztig-canonicalII}
\bysame, \emph{Canonical bases arising from quantized enveloping algebras.
  {II}}, no. 102, 1990, Common trends in mathematics and quantum field theories
  (Kyoto, 1990), pp.~175--201.

\bibitem[Lus10]{lusztig-book}
\bysame, \emph{Introduction to quantum groups}, Modern Birkh\"{a}user Classics,
  Birkh\"{a}user/Springer, New York, 2010, Reprint of the 1994 edition.

\bibitem[Mac03]{macdonald-affine}
I.~G. Macdonald, \emph{Affine {H}ecke algebras and orthogonal polynomials},
  Cambridge Tracts in Mathematics, vol. 157, Cambridge University Press,
  Cambridge, 2003.

\bibitem[Nak01]{nakajima-qaffine}
H.~Nakajima, \emph{Quiver varieties and finite-dimensional representations of
  quantum affine algebras}, J. Amer. Math. Soc. \textbf{14} (2001), no.~1,
  145--238.

\bibitem[Sai94]{saito-PBW}
Y.~Saito, \emph{P{BW} basis of quantized universal enveloping algebras}, Publ.
  Res. Inst. Math. Sci. \textbf{30} (1994), no.~2, 209--232.

\bibitem[So{\u{\i}}90]{soibelman}
Ya.~S. So{\u{\i}}belman, \emph{Algebra of functions on a compact quantum group
  and its representations}, Algebra i Analiz \textbf{2} (1990), no.~1,
  190--212.

\bibitem[TL11]{valerio3}
V.~Toledano~Laredo, \emph{The trigonometric {C}asimir connection of a simple
  {L}ie algebra}, J. Algebra \textbf{329} (2011), 286--327.


\end{thebibliography}

\newcommand{\etalchar}[1]{$^{#1}$}
\providecommand{\bysame}{\leavevmode\hbox to3em{\hrulefill}\thinspace}
\providecommand{\MR}{\relax\ifhmode\unskip\space\fi MR }
\providecommand{\MRhref}[2]{%
  \href{http://www.ams.org/mathscinet-getitem?mr=#1}{#2}
}
\providecommand{\href}[2]{#2}

\end{document}